\documentclass{m2an}

\usepackage{amsmath,amssymb}
\usepackage[latin1]{inputenc}
\usepackage{epsfig,float,psfrag}
\usepackage{color}
\newtheorem{thm}{Theorem}

\newtheorem{alem}[thm]{Lemma}
\newtheorem{aprop}[thm]{Proposition}
\newtheorem{acor}[thm]{Corollary}
\newtheorem{arem}[thm]{Remark}
\newenvironment{adem}[1][]%
   {\ \\ {\bf Proof #1~: }}%
   {\hfill\mbox{\rule{2 true mm}{3 true mm}}\vskip 2 ex\noindent}
\newcommand{\E}{{\mathbb E}}
\renewcommand{\P}{{\mathbb P}}
\newcommand{\R}{{\mathbb R}}
\newcommand{\N}{{\mathbb N}}
\renewcommand{\t}{{\theta}}
\newcommand{\dt}{{\delta t}}
\newcommand{\Dt}{{\Delta t}}
\newcommand{\uu}[1]  {{\boldsymbol{#1}} }
\newcommand{\ds}[1]{\displaystyle{#1}}

\begin{document}

\title{Diffusion Monte Carlo method : Numerical Analysis in a Simple
  Case}

\thanks{We thank Eric Cancès (CERMICS) and Mathias Rousset (Université
  Paul Sabatier, Toulouse) for many fruitful discussions and
  Michel Caffarel (IRSAMC, Toulouse) for introducing us to the DMC
  method and suggesting the toy model studied in this paper. We also
  thank the referees for useful suggestions which helped us to improve the first draft
  of this paper.}
\author{M.~El Makrini}\author{B.~Jourdain}\author{T.~Leli\`evre}\address{ENPC-CERMICS, 6-8 av Blaise Pascal, Cit\'e
    Descartes, Champs sur Marne, 77455 Marne-la-Vall\'ee Cedex 2, France
    -
    e-mail: \{makrimo,jourdain,lelievre\}@cermics.enpc.fr}


\begin{abstract}
The Diffusion Monte Carlo method is devoted to the computation of
electronic ground-state energies of molecules. In this paper, we focus on
implementations of this method which consist in exploring the
configuration space with a {\bf fixed} number of random walkers evolving
according to a Stochastic Differential Equation discretized in time. We
allow stochastic reconfigurations of the walkers to reduce the
discrepancy between the weights that they carry.  On a simple
one-dimensional example, we prove the convergence of the method for a
fixed number of reconfigurations when the number of walkers tends to
$+\infty$ while the timestep tends to $0$. We confirm our theoretical
rates of convergence by numerical experiments. Various resampling
algorithms are investigated, both theoretically and numerically.
\end{abstract}



\maketitle

\section*{Introduction}


The computation of electronic structures of atoms, molecules and solids is a
central problem in chemistry and physics. We focus here on electronic
ground state calculations where the objective is the computation of
 the lowest eigenvalue (the so-called ground-state
energy) $E_0$ of a self-adjoint Hamiltonian $H=-\frac{1}{2}\Delta +V$ with domain $D_\mathcal{H}(H)$ on a
Hilbert space $\mathcal{H} \subset L^2(\R^{3N})$ where $N$ is the number
of electrons (see~\cite{handbook} for a general introduction):
\begin{equation}\label{eq:E0}
E_0=\inf \{  \langle \psi,H\psi \rangle,\, \psi \in D_\mathcal{H}(H), \|\psi\|=1
\},
\end{equation}
where $\langle\cdot,\cdot\rangle$ denotes the duality bracket on $L^2(\R^{3N})$
and $\|\cdot\|$ the $L^2(\R^{3N})$-norm. For simplicity, we omit the spin variables. The function $V$ describes the interaction between the electrons, and between the
electrons and the nuclei, which are supposed to be fixed point-like particles. The functions $\psi$ are square
integrable, their normalized square modulus $|\psi|^2$ being interpreted as the
probability density of the particles positions in space, and they satisfy an antisymmetry condition with respect to the
numbering of the electrons, due to the fermionic nature of the
electrons (Pauli principle): $\mathcal{H}=\bigwedge_{i=1}^N L^2(\R^3)$. We suppose that the potential $V$ is such that $E_0$ is an isolated
eigenvalue of $H$ (see~\cite{CJL} for sufficient conditions), and we
denote by $\psi_0$ a normalized eigenfunction
associated with $E_0$. 


Due to the
high dimensionality of the problem, stochastic methods are particularly
well suited to compute~$E_0$. The first particle approximation
scheme of such spectral quantities was introduced in~\cite{H84} for
finite state space models. Convergence analysis for such
interacting particle systems (both continuous
or discrete in time) first appeared in~\cite{DM03,Delmoral,Delmoralmiclo,DD04}. The Diffusion Monte Carlo (DMC) method is
widely used in chemistry (see~\cite{Caffarel,Umrigar}), but has been  only
recently considered from a mathematical viewpoint
(see~\cite{CJL,Rousset}). This method gives an estimate of~$E_0$ in
terms of the long-time limit of the expectation of a functional of
a drift-diffusion process with a source term. It requires an importance
sampling function~$\psi_I$ which approximates the ground-state $\psi_0$
of~$H$. Let us define the drift function $\uu{b}=\nabla \ln |\psi_I|$, the
so-called local energy $\ds{E_L=\frac{H \psi_I}{\psi_I}}$ and the DMC energy:
\begin{equation}\label{eq:EDMC}
E_{\rm DMC}(t)= \frac{\E\left(E_L(\uu{X}_t)\exp\left(-\int_0^t
      E_L(\uu{X}_s)ds\right)\right)}{\E\left(\exp\left(-\int_0^t
      E_L(\uu{X}_s)ds\right)\right)},
\end{equation}
where the $3N$-dimensional process $\uu{X}_t$ satisfies the stochastic
differential equation:
\begin{equation}\label{eq:SDE}
\left\lbrace
\begin{array}{l}
\ds{\uu{X}_t=\uu{X}_0 + \int_0^t \uu{b}(\uu{X}_s) \, ds + d \uu{W}_t,}\\
\uu{X}_0 \sim |\psi_I|^2(\uu{x}) \, d \uu{x}.
\end{array}
\right.
\end{equation}
The stochastic process $(\uu{W}_t)_{t \geq 0}$ is a standard $3N$-dimensional Brownian
motion. One can then show that (see~\cite{CJL})
\begin{equation}\label{eq:CV_E_DMC}
\lim_{t \to \infty} E_{\rm DMC}(t)=E_{{\rm DMC},0},
\end{equation}
where 
\begin{equation}\label{eq:EDMC0}
E_{{\rm DMC},0}=\inf \{  \langle
\psi,H\psi \rangle,\, \psi \in D_\mathcal{H}(H), \|\psi\|=1, \, \psi=0\mbox{ on }\psi_I^{-1}(0)
\}.
\end{equation}
We have proved in~\cite{CJL} that $E_{{\rm DMC},0} \geq E_0$, with equality if and
only if the nodal surfaces of $\psi_I$ coincide with those of a ground
state $\psi_0$ of $H$. In other words, if there exists a ground state
$\psi_0$ such that $\psi_I^{-1}(0)=\psi_0^{-1}(0)$, then $\lim_{t \to
  \infty} E_{\rm DMC}(t)=E_0$. The error $|E_0-E_{{\rm DMC},0}|$ is
related to the so-called fixed-node approximation, which is
well known by practitioners of the field (see~\cite{handbook}).

In this paper, we complement the theoretical results obtained in~\cite{CJL} with a numerical analysis in a simple case.
In practice, the longtime limit $E_{{\rm DMC},0}$ in (\ref{eq:CV_E_DMC}) is approximated by taking the
value of $E_{\rm DMC}$ at a (large) time $T>0$. Then $E_{\rm DMC}(T)$ is
approximated by using a
discretization in time of the stochastic differential
equation~(\ref{eq:SDE}) and of the integral in the exponential factor in
(\ref{eq:EDMC}), and an approximation of the expectation
values in~(\ref{eq:EDMC}) by an empirical mean over a large number $N$
of 
trajectories. These trajectories $(\uu{X}^i)_{1 \leq i \leq N}$, also called walkers in the
   physical literature or particles in the mathematical literature, satisfy a discretized version of~(\ref{eq:SDE}), and
   interact at times $n\Delta t$ for $n\in\{1,\hdots,\nu-1\}$ where
   $\Dt=T/\nu$ for $\nu\in{\mathbb N}^*$ through a stochastic
   reconfiguration step aimed at
   reducing the discrepancy between their exponential weights. We thus obtain an interacting
   particle system. The number of
   reconfiguration steps is $\nu-1$. The
stochastic differential equation \eqref{eq:SDE} is discretized
with a possibly smaller timestep $\dt=\Dt/\kappa=T/(\nu\kappa)$ with $\kappa\in\N^*$. The total
number of steps for the discretization of \eqref{eq:SDE} is then $K=\nu\kappa$.

In the following, we consider the following adapted version of the DMC
scheme with a fixed number of walkers (see~\cite{Caffarel}):
\begin{itemize}
\item[{$ \quad \bullet \;$}] {\bf Initialization} of an ensemble of $N$
  walkers $\left(\uu{X}^j_{0 \Dt} \right)_{1 \le j \le N}$ i.i.d. according to
  $|\psi_I|^2(\uu{x}) \, d\uu{x}$. 
\item[{$ \quad \bullet \;$}] {\bf Iterations in time:} let us be given
  the particle positions $\left( \uu{X}^j_{n \Dt} \right)_{1 \le j \le N}$
  at time $n \Delta t$, for $n \in \{0,\hdots,\nu-1\}$. The new particle positions at time $(n+1) \Delta t$ are obtained in two steps:
\begin{enumerate}
\item {\bf Walkers displacement:} for all $1 \le j \le N$, the
  successive positions $\left( \uu{X}^j_{n \Dt + \dt}, \, \hdots ,\,
    \uu{X}^j_{n \Dt + \kappa \dt}\right)$ over the time interval $(n \Delta t,(n+1) \Delta t)$ are obtained by an
appropriate discretization of~(\ref{eq:SDE}). In the field of
interacting particles system for Feynman-Kac formulae (see~\cite{Delmoral,Delmoralmiclo}), this step is called {\bf the mutation step}.
\item {\bf Stochastic reconfiguration:} The new positions\footnote{With a slight
  abuse of notation and though $n
    \Dt + \kappa \dt=(n+1)\Dt$, we distinguish between the particle positions $\uu{X}^j_{n
    \Dt + \kappa \dt}$ at the end of the walkers displacement on time interval
  $(n \Dt,(n+1) \Dt)$, and the new
  particle positions $\uu{X}^j_{(n+1)
    \Dt}$ obtained after the reconfiguration step, and
  which are used as the initial position for the next walkers displacement on time interval
  $((n+1) \Dt,(n+2) \Dt)$. We will use a more precise notation for the
  analysis of the numerical scheme in Section~\ref{sec:num_anal}, but this is not required
at this stage.} $\left(
    \uu{X}^j_{(n+1) \Dt} \right)_{1 \le j \le N}$ which will be used as
  the initial particle positions on the time interval $((n+1) \Dt,(n+2)
  \Dt)$ are obtained from independent sampling of the measure
\begin{equation}
\frac{\sum_{j=1}^N \exp\left( - \dt \sum_{k=1}^\kappa E_L(\uu{X}^j_{n
    \Dt + k \dt}) \right) \delta_{\uu{X}^j_{n
    \Dt + \kappa \dt}}}{ \sum_{j=1}^N \exp\left( - \dt \sum_{k=1}^\kappa E_L(\uu{X}^j_{n
    \Dt + k \dt}) \right)}.
\end{equation}
In words, the new particle positions $\left(
    \uu{X}^j_{(n+1) \Dt} \right)_{1 \le j \le N}$ are randomly chosen among the
  final particle positions $\left(
    \uu{X}^j_{n \Dt + \kappa \dt} \right)_{1 \le j \le N}$, each of them
  being weighted with the coefficient $\exp\left( - \dt \sum_{k=1}^\kappa E_L(\uu{X}^j_{n
    \Dt + k \dt}) \right)$ (accordingly to the exponential factor in~(\ref{eq:EDMC})).
 In the field of
interacting particles system for Feynman-Kac formulae, this step is called {\bf the
 selection step}.
\end{enumerate}
\end{itemize}
An estimate of $E_{\rm DMC}(t_{n+1})$ is then given by:
\begin{equation}\label{eq:EDMCapprox}
E_{\rm DMC}(t_{n+1}) \simeq \frac{1}{N} \sum_{j=1}^{N} \,
  E_L \left(\uu{X}^j_{(n+1) \Dt}\right).
\end{equation}
There are other possible estimations of $E_{\rm
  DMC}(t_{n+1})$. In~\cite{Caffarel}, the authors propose to use Cesaro
or weighted Cesaro means of the expression~(\ref{eq:EDMCapprox}). In Section~\ref{sec:num_anal}, we will use the
following expression: 
\begin{equation} \label{eq:EDMCapprox_prime}
E_{\rm DMC}(t_{n+1}) \simeq \frac{\sum_{j=1}^{N} \,
  E_L(\uu{X}^j_{n \Dt + \kappa \dt}) \exp\left( - \dt
    \sum_{k=1}^{\kappa} E_L(\uu{X}^j_{n \Dt + k \dt}) \right)}{\sum_{j=1}^{N}
  \exp\left( - \dt \sum_{k=1}^{\kappa} E_L(\uu{X}^j_{n \Dt + k \dt}) \right)},
\end{equation}
in an intermediate step to prove the convergence result. 

We would like to mention that a continuous in time version of the DMC scheme with stochastic
reconfiguration has been proposed in~\cite{Rousset}. The author
analyzes the longtime behavior of the interacting particle system and
proves in particular a uniform in time control of the variance of the
estimated energy.

The DMC algorithm presented above is prototypical. Many refinements
are used in practice. For example, an acception-rejection step is generally used in the
walkers displacement step (see~\cite{RCAL82}). This will not be
discussed here. Likewise, the selection
step can be done in many ways (see~\cite{CDM05,C04} for general
algorithms, and~\cite{Caffarel,Umrigar,Sorella} for algorithms used in
the context of DMC computations). In this paper, we restrict ourselves
to resampling methods with a fixed number of particles, and such that
the weights of the particles after resampling are equal to $1$. Then, the basic
consistency requirement of the selection step is that, conditionally on
the former positions $\left( \uu{X}^j_{n \Dt + k \dt} \right)_{1 \le
  j \le N, 1 \leq k \leq \kappa}$, the $i$-th particle~$\uu{X}^i_{n
    \Dt + \kappa \dt}$ is replicated $N \rho^i_n$ times in mean, where $\rho^i_n=\exp\left( - \dt \sum_{k=1}^\kappa E_L(\uu{X}^i_{n
    \Dt + k \dt}) \right) \Big/ \sum_{j=1}^N \exp\left( - \dt \sum_{k=1}^\kappa E_L(\uu{X}^j_{n
    \Dt + k \dt}) \right)$ denotes the (normalized) weight of the $i$-th particle. There are of course many ways to satisfy this requirement.
We presented above the so-called
{\it multinomial resampling} method. We will also discuss below {\it residual
resampling} (also called stochastic remainder resampling), {\it
stratified resampling} and {\it
systematic resampling}, which may also be used for DMC computations. Let us briefly describe these three resampling
methods. Residual resampling consists in reproducing
$\lfloor N \rho^i_n  \rfloor$ times the $i$-th particle, and then completing
the set of particles by using multinomial resampling to draw the
$N^R=N-\sum_{l=1}^N \lfloor N \rho^l_n  \rfloor$ remaining particles, the
$i$-th particle being assigned the weight $\rho^{R,i}_n=\{ N \rho^i_n \} /
N^R$. Here and in the following,
$\lfloor x\rfloor$ and $\{x\}$ respectively denote the
  integer and the fractional part of $x\in\R$. In the stratified
  resampling method, the interval $(0,1)$ is divided into $N$ intervals $((i-1)/N,i/N)$
  ($1 \leq i \leq N$), $N$ random variables are then drawn independently and
  uniformly in each interval, and the new particle positions are then
  obtained by the inversion method: $\uu{X}^i_{(n+1)\Dt} = \sum_{j=1}^N
  1_{\{\sum_{l=1}^{j-1} \rho^l_n < (i-U^i_n)/N \leq \sum_{l=1}^{j}
    \rho^l_n \}} \uu{X}^j_{n\dt + \kappa \dt} $, where $U^i_n$ are
  i.i.d. random variables uniformly distributed over $[0,1]$. Here and
  in the following, we use the convention $\sum_{l=1}^0 \cdot =0$.
Systematic resampling
  consists in replicating the $i$-th particle $\Big\lfloor N \sum_{l=1}^i
  \rho^l_n + U_n \Big\rfloor - \Big\lfloor N \sum_{l=1}^{i-1} \rho^l_n + U_n \Big\rfloor$
  times\footnote{The consistency of this resampling method follows
    from the following easy computation $$\E\left(\lfloor
      x+U\rfloor\right)=\lfloor x\rfloor\P(U<1-\{x\})+(\lfloor
  x\rfloor+1)\P(U\geq 1-\{x\})=\lfloor x\rfloor(1-\{x\})+(\lfloor
  x\rfloor+1)\{x\}=x.$$},  where $(U_n)_{n \geq 1}$ are independent random variables uniformly distributed in
  $[0,1]$. Notice that systematic resampling can be seen as the
  stratified resampling method, with $U^1_n= \ldots=U^N_n=U_n$.
Contrary to the three other resampling methods, after a systematic
resampling step,
 the new particle positions are not independent,
  conditionally on the former positions. This makes
  systematic resampling much
  more difficult to study mathematically. To our knowledge, its
  convergence even in a discrete time setting is still an open question. We will therefore restrict
  ourselves to a numerical study of its performance.

Notice that practitioners often use branching algorithms with an evolving number
of walkers during the computation (see~\cite{RCAL82,Umrigar}): the particles with low local
energy are replicated and the particles with high local energy are
killed, without keeping the total number of particles constant. This may
lead to a smaller Monte Carlo error (fourth contribution to the
error in the classification just below).



We can distinguish between four sources of errors in the approximation
of $E_0$ by $\ds{ \frac{1}{N} \sum_{j=1}^{N} \,
  E_L \left(\uu{X}^j_{\nu \Dt}\right)}$:
\begin{enumerate}
\item the error due to the fixed node approximation $|E_0-E_{{\rm
    DMC},0}|$,
\item the error due to finite time approximation of the limit: $\lim_{t
    \to \infty} E_{\rm DMC}(t) \simeq E_{\rm DMC}(T)$,
\item the error due to the time discretization of the stochastic differential
equation~(\ref{eq:SDE}) and of the integral in the exponential factor in
$E_{\rm DMC}(t)$ (see~(\ref{eq:EDMC})),
\item the error introduced by the interacting particle system, due to
  the approximation of the expectation value in~(\ref{eq:EDMC}) by an empirical mean.
\end{enumerate}

The error (1) due to the fixed node approximation has been analyzed
theoretically in~\cite{CJL}. 

Concerning the error (2) due to finite time approximation of the limit, the rate of convergence in time is
typically exponential. Indeed if $H$ admits a spectral gap (namely if the
distance between $E_0$ and the remaining of the spectrum of $H$ is strictly
positive), and if $\psi_I$ is such that $\langle \psi_I, H \psi_I \rangle <
\inf \sigma_{\rm ess}(H)$, then one can show that the operator $H$ with
domain $D_\mathcal{H}(H) \cap \{\psi, \, \psi=0\mbox{ on }\psi_I^{-1}(0) \}$ (whose
lowest eigenvalue is $E_{{\rm DMC},0}$, see~(\ref{eq:EDMC0})) also admits a
spectral gap $\gamma>0$. Then, by standard spectral decomposition methods,
we have:
$$0 \leq \left| E_{\rm DMC}(t) - E_{{\rm DMC},0} \right| \leq C \exp(-\gamma t).$$

Our aim in this paper is to provide some theoretical and numerical
results related to the errors~(3) and~(4), in the framework
of a simple one-dimensional case. We therefore consider in the
following that the final time of simulation $T$ is fixed and we analyze the error introduced by the
numerical scheme on the estimate of $E_{\rm DMC}(T)$. Our convergence
result is of the form:
\begin{equation}\label{eq:est_err_general}
\E \left|E_{\rm DMC}(T)- \frac{1}{N} \sum_{j=1}^{N} \,
  E_L \left(\uu{X}^j_{\nu \kappa \dt}\right)\right| \leq C(T) \, \dt +
\frac{C(T,\nu)}{\sqrt{N}},
\end{equation}
where $C(T)$ (resp. $C(T,\nu)$) denotes a constant which only depends on
$T$ (resp. on $T$ and $\nu$) (see Theorem~\ref{th:est_err} and
Corollary~\ref{cor:est_err} below).


Let us now present the toy model we consider in the following. We consider the Hamiltonian 
\begin{equation}\label{eq:H}
H=-\frac{1}{2}\frac{d^2}{dx^2}+V, \mbox{ with }
V=\frac{\omega^2}{2}x^2+\theta x^4,
\end{equation}
where $\omega,\theta>0$ are two constants. The ground
state energy $E_0$ is defined by~(\ref{eq:E0}), with
\begin{equation}\label{eq:cal_H}
\mathcal{H}=\left\{ \psi \in L^2(\R), \, \psi(x)=-\psi(-x)
\right\}.
\end{equation}
We restrict the functional spaces to odd functions in order to mimic the antisymmetry
constraint on $\psi$ for fermionic systems. The importance sampling
$\psi_I$ is chosen to be the ground state of
$H_0=-\frac{1}{2}\frac{d^2}{dx^2} + \frac{\omega^2}{2}x^2$ on
$\mathcal{H}$:
\begin{equation}\label{eq:psi_I}
\psi_I(x)=\sqrt{2\omega}\left(\frac{\omega}{\pi}\right)^{1/4}xe^{-\frac{\omega
    }{2}x^2}.
\end{equation}
It is associated with the energy $\frac{3}{2}\omega$: $H_0 \psi_I=
\frac{3}{2}\omega \psi_I$.
The drift function $b$ and the local energy $E_L$ are then defined by:
\begin{equation}
b(x)=\frac{\psi_I'}{\psi_I}(x)=\frac{1}{x}-\omega x, \mbox{ and } E_L(x)=V(x)-\frac{1}{2}\;\frac{\psi_I''}{\psi_I}(x)=\frac{3}{2}\omega+\theta x^4.
\end{equation}
Thus, using equation~(\ref{eq:EDMC}), the DMC energy is:
\begin{equation}E_{\rm DMC}(t)=\frac{3}{2}\omega+\theta\frac{\E\left(X^4_t\exp\left(-\theta\int_0^t
      X_s^4ds\right)\right)}{\E\left(\exp\left(-\theta\int_0^t
      X_s^4ds\right)\right)},
\label{defedmc}\end{equation}
where
\begin{equation}
X_t=X_0 + \int_0^t \left(\frac{1}{X_s}-\omega X_s\right)ds+ W_t,\label{eds}
\end{equation}
with $(W_t)_{t \geq 0}$ a Brownian motion independent from the initial variable $X_0$
which is distributed according to the invariant measure $2\psi_I^2(x)1_{\{x>0\}}dx$.
We recall that due to the explosive part in the drift function~$b$, the stochastic
process cannot cross $0$, which is the zero point of $\psi_I$
(see~\cite{CJL}): $\P(\exists t>0,\, X_t=0)=0$. This explains why the restriction of $\psi_I^2$ to
$\R_+^*$ is indeed an invariant measure for~(\ref{eds}). For $\theta>0$,
the longtime
limit $E_{{\rm DMC},0}$ of $E_{\rm DMC}(t)$ is not analytically known,
but can be very accurately computed by a spectral method (see
Section~\ref{sec:num_res_spec}). Let us finally make precise that for the
numerical analysis, we use a special feature of our simple model, namely
the fact that for $s \leq t$, it is possible to simulate the conditional
law of $X_t$ given $X_s$ (see Appendix). The time discretization error is thus only related to
the discretization of the integral in the exponential factor in the DMC
energy~(\ref{eq:EDMC}). We however indicate some possible ways to
prove~(\ref{eq:E0}) with a convenient time discretization of the SDE (see Equation~(\ref{schemaurel}),
Remark~\ref{rem:weak} and Proposition~\ref{prop:Markov}).

 Though our model
is one-dimensional (and therefore still far from the real problem~(\ref{eq:E0})), it contains one of the main difficulties related to the
approximation of the ground state energy for fermionic systems, namely
the explosive behavior of the drift in the stochastic differential
equation. However, two characteristics of practical problems are missing
in the toy model considered here. First, since we consider a
one-particle model, we do not treat difficulties related to
singularities of the drift and of the local energy at points where two
particles (either two electrons or one electron and one nucleus) coincide. Second, the local energy $E_L$ generally explodes at the nodes
of the trial wave function, and this is not the case on the simple example
we study
since the trial wave function is closely related to the exact ground
state. For an adaptation of the DMC algorithm to take care of these
singularities, we refer to~\cite{Umrigar}.
Despite the simplicity of the model studied in this paper, we think that the
convergence results we obtain and the mathematical tools we use
are prototypical for generalization to more complicated systems.

Compared to previous mathematical analysis of convergence for
interacting particle systems with stochastic
reconfiguration~\cite{DM03,Delmoral,Delmoralmiclo,DD04,Rousset}, our
study concentrates on the limit $\dt \to 0$ and $N \to \infty$ for a
fixed time~$T$, and on
the {\em influence of the time discretization error} in the
estimate~(\ref{eq:est_err_general}), where the test function~$E_L$ is
{\em unbounded}. It is actually important in our analysis that this
unbounded function $E_L$ also appears in the weights of the particles,
since it allows for specific estimates (see Lemmas~\ref{majomom} and~\ref{decroiss} below).


The paper is organized as follows. In Section~\ref{sec:num_anal}, we
prove the convergence result, by adapting
the methods of~\cite{Delmoral,Delmoralmiclo} to analyze the dependence
of the error on $\dt$. We then check the optimality of this theoretical result
by numerical experiments in Section~\ref{sec:num_res}, where we also
analyze numerically the dependence of the results on
various numerical parameters, including the number $(\nu-1)$ of
reconfiguration steps. From these numerical experiments, we propose a
simple heuristic method to choose the optimal number of reconfiguration steps.

{\bf Notation:} For any set of
random variables $(Y_i)_{i \in I}$, we denote by $\sigma((Y_i)_{i \in
  I})$ the sigma-field generated by these random variables. The parameters $\omega$ and $\theta$ are fixed positive constants.
By convention, any sum from one to zero is equal to zero: $\sum_{k=1}^0
\cdot =0$. Likewise, the subset $\{1,2, \ldots ,0\}$ of $\N$ is by convention the
empty set. For any real $x$, $\lfloor x\rfloor$ and $\{x\}$ respectively denote the
  integer and the fractional part of $x$.
\section{Numerical Analysis in a Simple Case}\label{sec:num_anal}

We perform the numerical analysis in two steps: time discretization and
then particle approximation.

\subsection{Time discretization}\label{sec:time_disc}

We recall that $T>0$ denotes the final simulation time, and that
$\dt=\frac{T}{K}$ is the smallest time-step. Since $Y_t=X_t^2$ is a
square root process solving $dY_t=(3-2\omega Y_t)dt+2\sqrt{Y_t}dW_t$,
it is possible to simulate the increments $Y_{(k+1)\dt}-Y_{k\dt}$ and
therefore 
$X_{(k+1)\dt}-X_{k\dt}$ (see Appendix or \cite{Glasserman} p.120). We
can thus simulate exactly in law the vector
$(X_0,X_{\dt},\hdots,X_{K \dt})$. That is why we are first going to study the
error related to the time discretization of the integral which appears
in the exponential factors in \eqref{defedmc}.\par
Let us define the corresponding approximation of $E_{\rm
   DMC}(T)$:
\begin{equation}\label{eq:EDMCdt}
E_{\rm
   DMC}^{\dt}(T)=\frac{\E\left(E_L(X_{T})\exp\left(-\dt\sum_{k=1}^KE_L(X_{k\dt})\right)\right)}{\E\left(\exp\left(-\dt\sum_{k=1}^KE_L(X_{k\dt})\right)\right)}=\frac{3}{2}\omega+\theta\frac{\E\left(X^4_T\exp\left(-\t\dt\sum_{k=1}^KX^4_{k\dt})\right)\right)}{\E\left(\exp\left(-\t\dt\sum_{k=1}^KX^4_{k\dt}\right)\right)}.
\end{equation}
\begin{aprop}
$$
  \forall K\in\N^*,\left|E_{\rm
   DMC}(T)-E_{\rm
   DMC}^{\dt}(T)\right|\leq C_T\dt.
$$
\label{erreurdt}\end{aprop}
\begin{adem}
Using H\"older inequality, we have:
\begin{align*}
   \left|E_{\rm DMC}(T)-E_{\rm DMC}^{\dt}(T)\right|\leq&\frac{\t}{\E\left(\exp\left(-\t\dt\sum_{k=1}^KX^4_{k\dt}\right)\right)}\left(\sqrt{\E(X_T^8)}+\frac{\E\left(X^4_{T}\exp\left(-\theta\int_0^T
      X_s^4ds\right)\right)}{\E\left(\exp\left(-\theta\int_0^T
      X_s^4ds\right)\right)}\right)\\
&\left(\E\left(\left(\exp\left(-\t\int_0^T
      X_s^4ds\right)-\exp\left(-\t\dt\sum_{k=1}^KX^4_{k\dt}\right)\right)^2\right)\right)^{1/2}.
\end{align*}
The conclusion is now a consequence of Lemma~\ref{vitfort} and the fact that the function $x\in\R_+\rightarrow e^{-\t x}$ is Lipschitz
continuous with constant $\t$.
\end{adem}
      
\begin{alem}
For any $K\in\N^*$,
$$\E\left(\left(\int_0^T
      X_s^4ds-\dt\sum_{k=1}^KX^4_{k\dt}\right)^2\right)\leq
  C\dt^2(T^2+T),$$
where $\dt=\frac{T}{K}$.
\label{vitfort}\end{alem}

\begin{adem}[of Lemma \ref{vitfort}]
   By Itô's formula, $dX^4_t=(10X^2_t-4\omega X^4_t)dt+4X_t^3dW_t$. With
   the integration by parts formula, one deduces that for any $k\in\{1,\hdots,K\}$,
$$\int_{(k-1)\dt}^{k\dt}(X^4_{k\dt}-X^4_s)ds=\int_{(k-1)\dt}^{k\dt}(s-(k-1)\dt)\left((10X^2_s-4\omega
  X^4_s)ds+4X_s^3dW_s\right).$$
Therefore denoting $\tau_s=\lfloor\frac{s}{\dt}\rfloor\dt$ the
  discretization time just before $s$, one obtains
$$\dt\sum_{k=1}^KX^4_{k\dt}-\int_0^T
      X_s^4ds=\int_0^T(s-\tau_s)(10X^2_s-4\omega
  X^4_s)ds+\int_0^T(s-\tau_s)4X_s^3dW_s.$$
Hence 
\begin{align*}
   \E\left(\left(\dt\sum_{k=1}^KX^4_{k\dt}-\int_0^T
      X_s^4ds\right)^2\right)\leq 2\int_0^T(s-\tau_s)^2\E\left(T(10X^2_s-4\omega
  X^4_s)^2+16X_s^6)\right)ds.
\end{align*}
Since $X_0$ is distributed according to the invariant measure
$2\psi_I^2(x)1_{\{x>0\}}dx$, so is $X_s$. As a consequence, for any $p\in\N$, $\E(X_s^p)$ does not
depend on $s$ and is finite and the conclusion follows readily.
\end{adem} 

In realistic situations, exact simulation of the increments
$X_{(k+1)\dt}-X_{k\dt}$ is not possible and one has to resort to
discretization schemes. The singularity of the drift coefficient
prevents the process $X_t$ from crossing the nodal surfaces of the
importance sampling function $\psi_I$. The standard explicit Euler
scheme does not preserve this property at the discretized level. For
that purpose, we suggest to use the following explicit scheme proposed by
\cite{Alfonsi}
\begin{equation}
   \begin{cases}
      \bar{X}_0=X_0,\\
\ds{\forall k\in\N,\;\bar{X}_{(k+1)\delta t}=\left(\left({\bar{X}_{k\delta t}}(1-\omega
    \delta t)+\frac{\Delta W_{k+1}}{1-\omega
    \delta t}\right)^2+2\delta t\right)^{1/2}}\;\;\mbox{with}\;\;\Delta
    W_{k+1}=W_{(k+1)\delta t}-W_{k\delta t}.
   \end{cases}\label{schemaurel}
\end{equation}
Because of the singularity at the origin of the drift coefficient in
\eqref{eds}, we have not been able so far to prove the following weak error
bound (see Remark~\ref{rem:weak} below):
\begin{equation}
   \left|E\left(f(X^4_{T})\exp\left(-\theta\int_0^{T}
      X_s^4ds\right)\right)-\E\left(f(\bar{X}^4_{T})\exp\left(-\theta \delta t\sum_{k=1}^{K}
      \bar{X}_{k\delta t}^4\right)\right)\right|\leq C_T\dt\mbox{ for
      }f(x)\equiv 1\mbox{ and }x^4.\label{weakaurel}
\end{equation}
Such a bound is expected according to~\cite{Talay-Tubaro} and would imply that
\begin{equation}
   \left|E_{\rm DMC}(T)-\frac{\E\left(E_L(\bar{X}_{T})\exp\left(-\dt\sum_{k=1}^KE_L(\bar{X}_{k\dt})\right)\right)}{\E\left(\exp\left(-\dt\sum_{k=1}^KE_L(\bar{X}_{k\dt})\right)\right)}\right|\leq C_T\dt.\label{discretaurel}
\end{equation}
\begin{arem}\label{rem:weak}
We would like to sketch a possible way to
prove~\eqref{weakaurel}. Because the square root in \eqref{schemaurel}
makes expansions with respect to $\delta t$ and $\Delta W_{k+1}$
complicated, it is easier to work with $Y_t=X_t^2$ and $\bar{Y}_{k\delta t}=\bar{X}^2_{k\delta
    t}$ which satisfy
$$dY_t=(3-2\omega Y_t)dt+2\sqrt{Y_t}\;dW_t\;\;\mbox{and}\;\;\bar{Y}_{(k+1)\delta t}=\left(\sqrt{\bar{Y}_{k\delta t}}(1-\omega
    \delta t)+\frac{\Delta W_{k+1}}{1-\omega
    \delta t}\right)^2+2\delta t.$$
The standard approach to analyze the time discretization error of the
numerator and denominator of the left hand side of~\eqref{discretaurel} is then to introduce some
functions $v$ and $w$ solutions to the partial differential equation:
\begin{equation}\label{eq:EDP}
\partial_t v=(3-2y)\partial_y
v+2y\partial_{yy}v-\t y^2v,\;(t,y)\in\R_+\times (0,+\infty)
\end{equation}
with initial conditions $v(0,y)=y^2$ and $w(0,y)=1$. Now, we write (for
the numerator, for example):
\begin{align*}
&\E\left(X^4_{T}\exp\left(-\theta\int_0^{T}
      X_s^4ds\right)\right)-\E\left(\bar{X}^4_{T}\exp\left(-\theta \delta t\sum_{k=1}^{K}
      \bar{X}_{k\delta t}^4\right)\right)\\
&=\sum_{k=0}^{K-1}\E\left(\left(v(T-k\delta t,\bar{Y}_{k\delta
      t})-e^{-\theta\delta t\bar{Y}_{(k+1)\delta t}^2}v(T-(k+1)\delta t,\bar{Y}_{(k+1)\delta
      t})\right)\exp\left(-\theta \delta t\sum_{j=0}^{k-1}
      \bar{Y}_{j\delta t}^2\right)\right).
\end{align*}
An error bound of the form $C_T \dt$ can now be proved by some Taylor
expansions as in \cite{Talay-Tubaro,Alfonsi}, provided
the existence of a sufficiently smooth solution $v$
to~(\ref{eq:EDP}). We have not been able to prove existence of such a
solution so far.
\end{arem}

\subsection{Particle approximation}

\begin{figure}[htbp]
\input{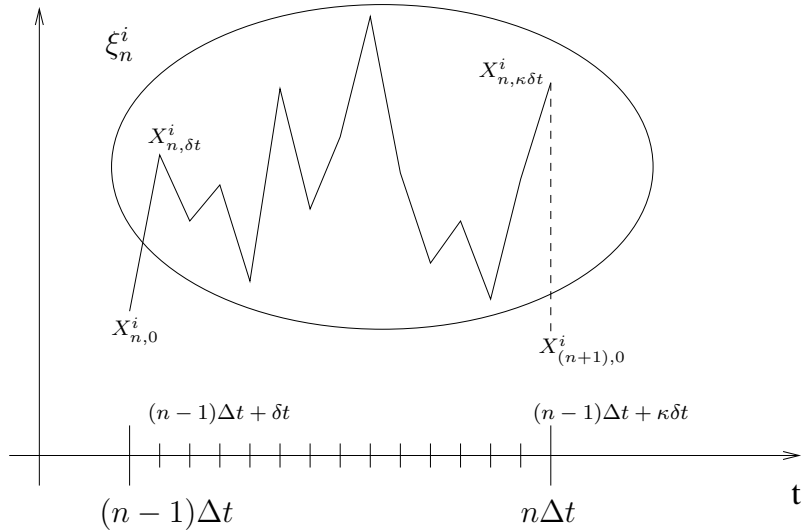}
\caption{The $i$-th particle $\xi^i_n$ at time index $n$ is composed of the
  successive positions $(X^i_{n,\dt},\hdots,X^i_{n,\kappa \dt})$ of the
  $i$-th walker on time interval $((n-1) \Dt,n \Dt)$.}\label{fig:particle}
\end{figure}

We now introduce some notation to study the particle approximation. We
recall that $\nu$ denotes the number of large timesteps (the
number of reconfiguration steps is $\nu-1$), and
$\Dt=\kappa \dt$ the time period between two reconfiguration steps. Let
us suppose that we know the initial positions $(X^i_{n,0})_{1 \leq i \leq N}$ of
the $N$ walkers at time $(n-1)\Dt$, for a time index $n \in \{1,\hdots,\nu\}$. The successive positions of the walkers
over the time interval $((n-1) \Dt,n\Dt)$ are then given by
$(X^i_{n,\dt},\hdots,X^i_{n,\kappa\dt})$, where $(X^i_{n,t})_{0 \leq t \leq
\Dt}$ satisfies:
\begin{equation}
X^i_{n,t}=X^i_{n,0}+\int_0^t b(X^i_{n,s}) \, ds+ \left(W^i_{t+(n-1) \Dt} -
  W^i_{(n-1) \Dt}\right).\label{sden+1}
\end{equation}
Here $(W^1,\hdots,W^N)$ denotes a $N$-dimensional Brownian motion independent
from the initial positions of the walkers $(X^i_{1,0})_{1\leq i\leq N}$
which are i.i.d. according
to $2 \psi_I^2(x) 1_{\{x>0\}} dx$. We recall that in our framework, it
is possible to simulate exactly in law all these random variables (see
Appendix). We store the successive positions
$(X^i_{n,\dt},\hdots,X^i_{n,\kappa\dt})$ of the $i$-th walker over the
time interval $((n-1) \Dt,n \Dt)$ in a so-called particle $\xi^i_{n}
\in (\R_+^*)^\kappa$
(see Figure~\ref{fig:particle}): $\forall i \in \{1,\hdots,N\}, \forall  n
\in \{1,\hdots, \nu\}$, 
\begin{equation}
\xi^i_{n}=(X^i_{n,\dt},\hdots,X^i_{n,\kappa\dt}).
\end{equation}
In the following, we will denote by $\xi_{n}=(\xi^1_{n},\hdots,\xi^N_{n})$ the configuration of the
ensemble of particles at time index $n$. We have here described {\bf the mutation step}.

For a given configuration of the particles $\xi_{n}$ at a time index $n \in \{1,\hdots, \nu\}$, {\bf the selection step}
now consists in choosing the initial positions
$(X^i_{n+1,0})_{1 \leq i \leq N}$ of the $N$ walkers at time $n \Dt$
using one of the following resampling algorithm:
\begin{itemize}
   \item[{\bf(S1)}]
The $(X^i_{n+1,0})_{1 \leq i \leq N}$ are conditionally
independent w.r.t.~$\xi_{n}$ and for $1 \leq i \leq N$, $X^i_{n+1,0}$ is distributed according to the measure
 \begin{equation}
   \epsilon_n
   g(\xi_n^i)\delta_{\xi_{n,\kappa}^i}+(1-\epsilon_n
   g(\xi_n^i)) \sum_{j=1}^N \rho_n^j \delta_{\xi_{n,\kappa}^j},\label{messel}
\end{equation}
where $g$ is defined by, for $y=(y_1,\hdots,y_\kappa)\in(\R_+^*)^\kappa$,
\begin{equation}\label{eq:g}
g(y)=\exp\left(-\t\dt\sum_{k=1}^\kappa y_k^4\right),
\end{equation}
$\rho_n^j$ denotes the weight of the $j$-th particle
\begin{equation}\label{eq:rho}
\rho_n^j=\frac{g(\xi_n^j)}{\sum_{j=1}^Ng(\xi_n^j)}
\end{equation}
and $\epsilon_n$ is a non negative function of $\xi_{n}$ such that
$\epsilon_n \leq 1 \big/ \max_{1\leq
    i\leq N}g(\xi_n^i)$. In particular the following choices are possible for $\epsilon_n$:
\begin{equation}
   \epsilon_n=0,\;\epsilon_n=1\;\mbox{ or }\;\epsilon_n=\frac{1}{\max_{1\leq
    i\leq N}g(\xi_n^i)}.\label{choix_eps}
\end{equation}
The so-called {\bf multinomial resampling} method which corresponds to the
choice $\epsilon_n=0$ gives rise to a maximum decorrelation with the former
position of the particles, while with growing
$\epsilon_n$, more and more correlation is introduced.
\item[{\bf (S2)}] The $(X^i_{n+1,0})_{1 \leq i \leq N}$ are such that
\begin{equation}
   \begin{cases}
      \forall j\in\{1,\hdots,N\},\;\forall
i\in\left\{ \left(1+ \sum_{l=1}^{j-1} a^l_n\right),\hdots,\left(\sum_{l=1}^{j} a^l_n\right)\right\}, \\
\qquad X^i_{n+1,0}=\xi^j_{n,\kappa},\\
\mbox{and the variables $(X^i_{n+1,0})_{1+ \sum_{l=1}^{N} a^l_n \leq i
    \leq N}$ are conditionally
independent w.r.t.~$\xi_{n}$,}\\
\qquad \mbox{with } X^i_{n+1,0}\mbox{ distributed according to
}\sum_{j=1}^N\left\{ N \rho^j_n\right\}\delta_{\xi_{n,\kappa}^j}
\Big/ \left(N- \sum_{l=1}^{N} a^l_n\right),
   \end{cases}
\label{messel2}\end{equation}
where
  \begin{equation}
   a^j_n=\big\lfloor N \rho_n^j \big\rfloor,\;j\in\{1,\hdots,N\}.\label{defpoidsa}
  \end{equation} 
Notice that the $(X^i_{n+1,0})_{1 \leq i \leq N}$ are conditionally
independent w.r.t.~$\xi_{n}$.
This is the so-called {\bf residual resampling} method.
\item[{\bf (S3)}] The $(X^i_{n+1,0})_{1 \leq i \leq N}$ are such that, for $1 \leq i \leq N$,
\begin{equation}\label{eq:strat}
X^i_{(n+1),0} = \sum_{j=1}^N
  1_{\{\sum_{l=1}^{j-1} \rho^l_n < (i-U^i_n)/N \leq \sum_{l=1}^{j}
    \rho^l_n \}} \xi^j_{n,\kappa},
\end{equation}
where $(U^i_n)_{1 \leq i \leq N}$ are random variables
i.i.d. according to the uniform law on $[0,1]$, independently of $\xi_{n}$. 
Notice that the $(X^i_{n+1,0})_{1 \leq i \leq N}$ are conditionally
independent w.r.t.~$\xi_{n}$. This is the so-called {\bf stratified resampling} method.
\end{itemize}

For $n\in \{1,\hdots,\nu\}$, let us denote by 
\begin{equation}\label{eq:eta_n_}
\eta_n^N=\frac{1}{N}\sum_{i=1}^N\delta_{\xi^i_n}
\end{equation} 
the particle approximation of the measure $\eta_n$ defined by: $\forall
f:(\R_+^*)^\kappa\rightarrow \R\mbox{ bounded }$,
\begin{equation}\label{eq:eta_n}
\eta_n(f)=\frac{\E\left(f\left(X_{(n-1)\Dt+\dt},\hdots,X_{(n-1)\Dt+\kappa\dt}\right)\exp\left(-\t\dt\sum_{k=1}^{(n-1)\kappa}(X_{k\dt})^4\right)\right)}{\E\left(\exp\left(-\t\dt\sum_{k=1}^{(n-1)\kappa}(X_{k\dt})^4\right)\right)},
\end{equation}
where the process $(X_t)_{0 \leq t \leq T}$ is defined by~\eqref{eds}.

For $y=(y_1,\hdots,y_\kappa)\in(\R_+^*)^\kappa$ and
$f:(\R_+^*)^\kappa \rightarrow \R$, we set 
\begin{equation}\label{eq:P}
Pf(y)=\E\left(f(X^{y_\kappa}_\dt,\hdots,X^{y_\kappa}_{\kappa \dt})\right)
\end{equation} 
where for $x\in\R^*_+$, 
\begin{equation}\label{eq:sdex}
X^x_t=x+\int_0^tb(X^x_s)ds+W_t
\end{equation}
 denotes the solution of
the stochastic differential equation~\eqref{eds} starting from~$x$.
By the Markov property, the measures $(\eta_n)_{1 \leq n \leq \nu}$
satisfy the inductive relations, for any function $f:(\R_+^*)^\kappa
\rightarrow \R$ bounded, $\forall n\in \{1,\hdots,\nu-1\}$,
\begin{align}
 \eta_{n+1}(f)&=\frac{\E\left(\exp\left(-\t\dt\sum_{k=1}^{n\kappa}(X_{k\dt})^4\right)\E\left(f\left(X_{n\Dt+\dt},\hdots,X_{n\Dt+\kappa\dt}\right)\bigg|(X_{j\dt})_{0\leq
 j\leq
 n\kappa}\right)\right)}{\eta_n(g)\E\left(\exp\left(-\t\dt\sum_{k=1}^{(n-1)\kappa}(X_{k\dt})^4\right)\right)}\\
&=\frac{1}{\eta_n(g)}\times\frac{\E\left(gPf\left(X_{(n-1)\Dt+\dt},\hdots,X_{(n-1)\Dt+\kappa\dt}\right)\exp\left(-\t\dt\sum_{k=1}^{(n-1)\kappa}(X_{k\dt})^4\right)\right)}{\E\left(\exp\left(-\t\dt\sum_{k=1}^{(n-1)\kappa}(X_{k\dt})^4\right)\right)}=\frac{\eta_n(gPf)}{\eta_n(g)},
   \label{inducetan}
\end{align}
where $g$ is defined by~(\ref{eq:g}). Moreover, we can express $E_{\rm
  DMC}^{\dt}(T)$ defined by~(\ref{eq:EDMCdt}) as:
\begin{equation}
   E_{\rm DMC}^{\dt}(T)=\frac{3}{2}\omega
+\theta\frac{\eta_\nu(gy_\kappa^4)}{\eta_\nu(g)}.\label{expr1}
\end{equation}
Therefore the particle approximation of $E_{\rm DMC}(T)$ is given by
\begin{equation}
   E_{\rm DMC}^{N,\nu,\kappa}(T)=\frac{3}{2}\omega
+\theta\frac{\eta^N_\nu(gy_\kappa^4)}{\eta^N_\nu(g)}.\label{expr2}
\end{equation}
This approximation of $E_{\rm DMC}(T)$ corresponds to the expression~(\ref{eq:EDMCapprox_prime})
given in the introduction. We will also prove in Corollary~\ref{cor:est_err} below the convergence of the
approximation which corresponds to the expression~(\ref{eq:EDMCapprox})
given in the introduction (see Equation~(\ref{expr3}) below).

The convergence of the approximation $E_{\rm DMC}^{N,\nu,\kappa}(T)$ is ensured by our main result :
\begin{thm}\label{th:est_err}
\begin{equation}\label{eq:est_err}
   \E\left|E_{\rm DMC}(T)-E_{\rm DMC}^{N,\nu,\kappa}(T)\right|\leq
   \frac{C}{\nu\kappa}+\frac{C_\nu}{\sqrt{N}},
\end{equation}
where the constant $C$ only depends on $T$ and the constant $C_\nu$ on $T$ and $\nu$.
\end{thm}
\begin{arem}\label{choinu}
The number of selection steps is $\nu-1$. For instance, when $\nu=1$,
   there is no selection involved in the expression of $E_{\rm
     DMC}^{N,\nu,\kappa}(T)$ and the particles remain independent. In
   this case, the first term in the right hand side
   of~(\ref{eq:est_err}) corresponds to the time discretization error
   proved in Proposition~\ref{erreurdt}, while the second term is
   the classical error estimate related to the law of large numbers. For a fixed number of
   selection steps, the theorem ensures the convergence of the particle
   approximation $E_{\rm DMC}^{N,\nu,\kappa}(T)$ as the time-step
   $\dt=T/(\nu \kappa)$ used for the discretization of the stochastic
   differential equation \eqref{eds} tends to~$0$ while the number $N$ of
   particles tends to $+\infty$. But this result does not specify the dependence of $C_\nu$ on $\nu$ and
   gives no hint on the optimal choice of the number of selection steps
   in terms of error minimization. We are going to deal with this
   important issue
   in the numerical study (see Section~\ref{sec:num_res}).
\end{arem}
According to the above expressions
   \eqref{expr1} and \eqref{expr2} of $E_{\rm DMC}^{\dt}(T)$ and $E_{\rm
   DMC}^{N,\nu,\kappa}(T)$, this theorem is easily proved by combining
   Proposition \ref{erreurdt} and the following result :
\begin{aprop}\label{propvitpart}
   \begin{equation}
   \E\left|\frac{\eta^N_\nu(gy_\kappa^4)}{\eta^N_\nu(g)}-\frac{\eta_\nu(gy_{\kappa}^4)}{\eta_\nu(g)}\right|\leq \frac{C_\nu}{\sqrt{N}}.\label{vitpart}
\end{equation}\end{aprop}

\begin{adem}[ of Proposition \ref{propvitpart}]
One has
\begin{align*}
   \E\left|\frac{\eta^N_\nu(gy_\kappa^4)}{\eta^N_\nu(g)}-\frac{\eta_\nu(gy_\kappa^4)}{\eta_\nu(g)}\right|\leq &\frac{\E|\eta_\nu^N(gy_\kappa^4)-\eta_\nu(gy_\kappa^4)|}{\eta_\nu(g)}\\&+\left(\E\left(\frac{\eta_\nu^N(gy_\kappa^4)}{\eta_\nu^N(g)}\right)^2\right)^{1/2}\frac{\left(\E\left(\eta_\nu^N(g)-\eta_\nu(g)\right)^2\right)^{1/2}}{\eta_\nu(g)}.
\end{align*}
According to Proposition \ref{vitmc} and Lemma \ref{contetag} below, the
first term of the right-hand-side and the quotient in the second term
are smaller than $C_\nu/\sqrt{N}$. Since by Jensen's inequality,
$\left(\frac{\eta_\nu^N(gy_\kappa^4)}{\eta_\nu^N(g)}\right)^2\leq
\frac{\eta_\nu^N(gy_\kappa^8)}{\eta_\nu^N(g)}$, the boundedness of 
$\E\left(\frac{\eta_\nu^N(gy_\kappa^4)}{\eta_\nu^N(g)}\right)^2$
follows from Lemma \ref{majomom} below. 
\end{adem}
\begin{aprop}\label{vitmc}
   For any bounded function $f:(\R_+^*)^\kappa\rightarrow \R$, 
\begin{equation}\label{eq:f_bounded}
\forall
n\in\{1,\hdots,\nu\},\;\E((\eta_n^N(f)-\eta_n(f))^2)\leq\frac{C_n}{N}\|f\|^2_\infty,
\end{equation}
where the constant $C_n$ does not depend on $\kappa$.

For any function $f:(\R_+^*)^\kappa\rightarrow \R$ such that for some $p\geq 2$,
$\ds{\|f\|_{\kappa,p}=\sup_{y\in\R_+^\kappa}\frac{|f(y)|}{1+y_\kappa^p}}$ is
finite,
\begin{equation}\label{eq:fkp_bounded}
 \forall
n\in\{1,\hdots,\nu\},\;\E|\eta_n^N(f)-\eta_n(f)|\leq\frac{C_n}{\sqrt{N}}\|f\|_{\kappa,p},
\end{equation}
where the constant $C_n$ does not depend on $\kappa$.
\end{aprop}
For $f$ bounded, the first estimate~(\ref{eq:f_bounded}) is proved in
\cite{Delmoralmiclo}. In order to prove Proposition~\ref{propvitpart}, we need to
apply Proposition~\ref{vitmc} with $f(y)=g(y)$ and $f(y)=g(y) y_\kappa^4$, which are
bounded functions with $L^\infty$ norm respectively equal to $1$ and
$\frac{C}{\dt}$ where $C$ is a constant not depending on $\dt$. But we
want to obtain the convergence when $\dt$ tends to $0$. This is why we need the
second estimate~(\ref{eq:fkp_bounded}), that we use with $f(y)=g(y)
y_\kappa^4$ for which $\|f\|_{\kappa,p}$ is bounded and does not depend on $\dt$.

Notice that for $f$ bounded, Corollary 2.20 in~\cite{Delmoralmiclo} states the convergence in law of
$\sqrt{N}(\eta^N_n(f)-\eta_n(f))$ to a centered Gaussian variable and
gives an expression of the variance of this limit variable. Because of
the complexity of this expression, using this result with $f(y)=g(y)y_\kappa^4$ did not really help us to
understand the dependence of $C_\nu$ on $\nu$ (see Remark \ref{choinu}
above).
\begin{adem}
For $f$ bounded, the first estimate~(\ref{eq:f_bounded}) is proved by
induction on $n$ in
\cite{Delmoralmiclo} (see Proposition 2.9). Since we follow the same
inductive reasoning to deal with $f$ such that
$\|f\|_{\kappa,p}<+\infty$, we give at the same time the proof for $f$ bounded.  

Since the initial positions $(\xi^i_1)_{1\leq i\leq N}$ are independent and identically
distributed with $\xi^i_{1,\kappa}$ distributed according to
$2\psi_I^2(x)1_{\{x>0\}}dx$, the
statement holds for $n=1$.\\
To deduce the statement at rank $n+1$ from the statement at rank $n$,
we remark that according to \eqref{inducetan}, 
\begin{equation}
   \eta_{n+1}^N(f)-\eta_{n+1}(f)=T_{n+1}+\frac{1}{\eta_n(g)}\left((\eta_n^N(gPf)-\eta_n(gPf))+\frac{\eta_n^N(gPf)}{\eta_n^N(g)}(\eta_n(g)-\eta_n^N(g))\right)\label{decomperr}
\end{equation}
where we recall that $P$ is defined by~(\ref{eq:P}), and 
$$T_{n+1}=\eta_{n+1}^N(f)-\frac{\eta_n^N(gPf)}{\eta_n^N(g)}.$$
To deal with this term $T_{n+1}$, one remarks that for the first type of
selection step (S1), all the possible choices
of $\epsilon_n$ given in \eqref{choix_eps} are $\sigma(\xi_n)$-measurable. As a consequence, for $i\in\{1,\hdots,N\}$,
$$\E(f(\xi^i_{n+1})|\xi_n)=\epsilon_n
g(\xi^i_n)Pf(\xi^i_n)+(1-\epsilon_n
g(\xi^i_n))\sum_{j=1}^N \rho^j_n Pf(\xi_{n}^j),$$
where $\rho^j_n$ is defined by~(\ref{eq:rho}). Multiplying
this equality by $\frac{1}{N}$ and summing over $i$, one deduces
\begin{equation}
   \E(\eta_{n+1}^N(f)|\xi_n)=\sum_{j=1}^N \rho_n^j  Pf (\xi^j_n)=\frac{\sum_{j=1}^Ng(\xi^j_n)Pf(\xi_n^j)}{\sum_{j=1}^N
   g(\xi^j_n)}=\frac{\eta_n^N(gPf)}{\eta_n^N(g)}.\label{condit}
\end{equation}
Now, for the stochastic remainder resampling algorithm (S2), by~\eqref{messel2},
$\E(\eta_{n+1}^N(f)|\xi_n)$ is equal to
$$\frac{1}{N}\sum_{j=1}^N\Bigg\lfloor\frac{Ng(\xi^j_n)}{\sum_{l=1}^N
  g(\xi^l_n)}\Bigg\rfloor
  Pf(\xi^j_n)+\sum_{i=1+\sum_{l=1}^N a_n^l}^N\frac{1}{N- \sum_{l=1}^N
    a_n^l} \, \sum_{j=1}^N\left\{\frac{Ng(\xi^j_n)}{\sum_{l=1}^N
  g(\xi^l_n)}\right\}Pf(\xi^j_n)$$
and \eqref{condit} still holds.  Finally, for the stratified resampling
method (S3), by~\eqref{eq:strat}, we have (using the footnote\footnotemark[2])
\begin{align*}
\E(\eta_{n+1}^N(f)|\xi_n)
&=\frac{1}{N}\sum_{i=1}^N\sum_{j=1}^N \E \left( 1_{\{\sum_{l=1}^{j-1} \rho^l_n < (i-U^i_n)/N \leq \sum_{l=1}^{j}
    \rho^l_n \}} \Big|\xi_n \right) Pf (\xi^j_n),\\
&=\frac{1}{N}\sum_{j=1}^N \E \left( \sum_{i=1}^N 1_{\{\sum_{l=1}^{j-1} \rho^l_n < (i-U^1_n)/N \leq \sum_{l=1}^{j}
    \rho^l_n \}} \Big|\xi_n \right) Pf (\xi^j_n),\\
&=\frac{1}{N}\sum_{j=1}^N \E \left( \left\lfloor   N \sum_{l=1}^{j}
    \rho^l_n  + U^1_n  \right\rfloor - \left\lfloor   N \sum_{l=1}^{j-1}
    \rho^l_n  + U^1_n  \right\rfloor \Big|\xi_n \right) Pf (\xi^j_n),\\
&=\sum_{j=1}^N \rho_n^j  Pf (\xi^j_n),
\end{align*}
which yields again~\eqref{condit}.
Since for all three possible selection steps, the variables $(\xi^i_{n+1})_{1\leq i\leq N}$ are independent
conditionally on $\xi_n$, one deduces that
\begin{align*}
   \E((T_{n+1})^2|\xi_n)=\frac{1}{N^2}\sum_{i=1}^N\E\left(\left(f(\xi_{n+1}^i)-\E(f(\xi^i_{n+1})|\xi_n)\right)^2|\xi_n\right)\leq \frac{1}{N}\E\left(\eta_{n+1}^N(f^2)|\xi_n\right).
\end{align*}
Therefore 
\begin{equation}
   \E((T_{n+1})^2)\leq \frac{1}{N}\E(\eta_{n+1}^N(f^2)).\label{majot}
\end{equation}
When $f$ is bounded, $\eta_{n+1}^N(f^2)\leq \|f\|_\infty^2$, $\left|\frac{\eta_n^N(gPf)}{\eta_n^N(g)}\right|\leq
\|Pf\|_\infty$, and $\|Pf\|_\infty\leq \|f\|_\infty$. Hence by \eqref{decomperr},
\begin{align*}
   \E((\eta_{n+1}^N(f)-\eta_{n+1}(f))^2)\leq 3\left(\frac{\|f\|_\infty^2}{N}+\frac{\E((\eta_n^N(gPf)-\eta_n(gPf))^2)+\|f\|_\infty^2\E((\eta_n^N(g)-\eta_n(g))^2)}{(\eta_n(g))^2}\right)
\end{align*}
with the second term of the right-hand-side smaller than $C\|f\|_\infty^2/N$ by the
induction hypothesis and Lemma \ref{contetag} below.\par
When $\|f\|_{\kappa,p}<+\infty$, combining \eqref{decomperr} and
\eqref{majot}, one obtains
\begin{align*}
   \E\left|\eta_{n+1}^N(f)-\eta_{n+1}(f)\right|\leq
   &\frac{\left(\E(\eta_{n+1}^N(f^2))\right)^{1/2}}{\sqrt{N}}+\frac{\E\left|\eta_n^N(gPf)-\eta_n(gPf)\right|}{\eta_n(g)}\\&+\left(\E\left(\frac{\eta_n^N(gPf)}{\eta_n^N(g)}\right)^2\right)^{1/2}\frac{ \left( \E(\eta_n^N(g)-\eta_n(g))^2 \right)^{1/2}}{\eta_n(g)}.
\end{align*}
Since $\|f^2\|_{k,2p} \leq 2 \|f\|_{k,p}^2$ (by using the inequality $f^2(y)\leq 2\|f\|_{\kappa,p}^2(1+y_\kappa^{2p})$), the first term of the right-hand-side is smaller than
$C_n\|f\|_{\kappa,p}/\sqrt{N}$ by Lemma \ref{majomom} below. Since, according to
Lemma \ref{contmom} below, $\|Pf\|_{\kappa,p}\leq e^{C_p\Delta t}\|f\|_{\kappa,p}$, the second term is smaller than
$C_n\|f\|_{\kappa,p}/\sqrt{N}$ by the induction hypothesis and Lemma
\ref{contetag}. Last, by using successively Cauchy Schwartz
inequalities, \eqref{condit} for $f^2$ and Lemma \ref{majomom}, one obtains that 
$\E\left(\frac{\eta_n^N(gPf)}{\eta_n^N(g)}\right)^2 \leq
\E\left(\frac{\eta_n^N(g (Pf)^2)}{\eta_n^N(g)}\right) \leq
\E\left(\frac{\eta_n^N(g P f^2)}{\eta_n^N(g)}\right) =
\E(\eta_{n+1}^N(f^2)) \leq C_n\|f\|^2_{\kappa,p}$. And it follows
from the Proposition statement for $f$ bounded and Lemma
\ref{contetag} that $\frac{\left(\E(\eta_n^N(g)-\eta_n(g))^2\right)^{1/2}}{\eta_n(g)}$ is smaller than
$C_n/\sqrt{N}$.
\end{adem}

\begin{arem}\label{rem:SRR}
Proposition~\ref{vitmc} (and therefore Theorem~\ref{th:est_err}) also
hold for the stratified remainder resampling algorithm, which consists
in combining the stochastic remainder resampling and the stratified
resampling. More precisely, it consists in replicating
$\lfloor N \rho^i_n  \rfloor$ times the $i$-th particle, and then completing
the set of particles by using stratified resampling to draw the
$N^R=N-\sum_{l=1}^N \lfloor N \rho^l_n  \rfloor$ remaining particles, the
$i$-th particle being assigned the weight $\rho^{R,i}_n=\{ N \rho^i_n \} /
N^R$.
\end{arem}

\begin{alem}\label{majomom} Let $h:(\R_+^*)^\kappa\rightarrow \R_+$ be such that for some
  $p\geq 2$, $\|h\|_{\kappa,p}<+\infty$. Then,
$$\forall n \in \{1,\hdots,\nu\},\; \max\left( \E(\eta_n^N(h)) ,
  \E\left(\frac{\eta_n^N(gh)}{\eta_n^N(g)}\right) \right)\leq
e^{C_pn\Dt}\|h\|_{\kappa,p}(1+\E(X_0)^p),$$
where $X_0$ is distributed according to the measure
$2\psi_I^2(x)1_{\{x>0\}}dx$ (see~(\ref{eds})).
\end{alem}
\begin{adem}
As the variables $\xi^i_{1,\kappa}, 1\leq i\leq N$ are distributed
according to the invariant measure $2\psi_I^2(x)1_{\{x>0\}}dx$, one has
$\E(\eta^N_1(h))\leq \|h\|_{\kappa,p}(1+\E(X_0)^p)$. In addition for $n\geq
1$, according to \eqref{condit},
$\E(\eta^N_{n+1}(h))=\E\left(\frac{\eta_n^N(gPh)}{\eta_n^N(g)}\right)$
where $\|Ph\|_{\kappa,p}\leq e^{C_p\Dt}\|h\|_{k,p}$  by Lemma \ref{contmom}. Therefore it is enough to
check the bound for $\E\left(\frac{\eta_n^N(gh)}{\eta_n^N(g)}\right)$.\par
For $n\geq 0$, one has
\begin{align}
 \E\left(\frac{\eta_{n+1}^N(gh)}{\eta_{n+1}^N(g)}\right) \leq
  \|h\|_{\kappa,p}\left(1+\E\left(\frac{\sum_{i=1}^N\exp\left(-\t\dt\sum_{k=1}^{\kappa }(\xi^i_{n+1,k})^4\right)(\xi^i_{n+1,\kappa})^p}{\sum_{j=1}^N\exp\left(-\t\dt\sum_{k=1}^{\kappa}(\xi^j_{n+1,k})^4\right)}\right)\right).\label{decompmom}
\end{align}
Let us denote in this proof $\xi^i_{n+1,0}=X^i_{n+1,0}$, where $0\leq n
\leq \nu-1$ and $1\leq i \leq N$. Let us set ${\mathcal
  F}=\sigma(\xi^i_{n+1,k},\;1\leq i\leq N,\;0\leq k\leq
  \kappa-1)$. By Lemma \ref{decroiss} below,
\begin{align*}
   \E\Bigg(&\frac{\sum_{i=1}^N\exp\left(-\t\dt\sum_{k=1}^{\kappa}(\xi^i_{n+1,k})^4\right)(\xi^i_{n+1,\kappa})^p}{\sum_{j=1}^N\exp\left(-\t\dt\sum_{k=1}^{\kappa}(\xi^j_{n+1,k})^4\right)}\bigg| {\mathcal
  F}\Bigg)\leq \frac{\sum_{i=1}^N\exp\left(-\t\dt\sum_{k=1}^{\kappa
   -1}(\xi^i_{n+1,k})^4\right)\E((\xi^i_{n+1,\kappa})^p|{\mathcal
  F})}{\sum_{j=1}^N\exp\left(-\t\dt\sum_{k=1}^{\kappa
   -1}(\xi^j_{n+1,k})^4\right)},\\
&\phantom{xxxxxxxxxxxxxxxxxxxxxx}=\frac{\sum_{i=1}^N\exp\left(-\t\dt\sum_{k=1}^{\kappa
   -1}(\xi^i_{n+1,k})^4\right)\E((X^x_\dt)^p)|_{x=\xi^i_{n+1,\kappa-1}}}{\sum_{j=1}^N\exp\left(-\t\dt\sum_{k=1}^{\kappa
   -1}(\xi^j_{n+1,k})^4\right)},\\
&\phantom{xxxxxxxxxxxxxxxxxxxxxx}\leq e^{C_p\dt}\frac{\sum_{i=1}^N\exp\left(-\t\dt\sum_{k=1}^{\kappa
   -1}(\xi^i_{n+1,k})^4\right)(\xi^i_{n+1,\kappa-1})^p}{\sum_{j=1}^N\exp\left(-\t\dt\sum_{k=1}^{\kappa
   -1}(\xi^j_{n+1,k})^4\right)}+e^{C_p\dt}-1,
\end{align*}
where we have used the definition of the mutation step (see~\eqref{sden+1}) and the Markov
   property for the stochastic differential equation \eqref{eq:sdex} to
   obtain the equality, and then Lemma \ref{contmom} for the last
   inequality. Notice that this estimate also holds for $\kappa=1$, in
   which case the right hand side reduces to $ \ds{\frac{e^{C_p\dt}}{N}
   (\xi^i_{n+1,0})^p +e^{C_p\dt}-1}$.

Taking expectations and iterating the reasoning, one deduces that
$$\E\left(\frac{\sum_{i=1}^N\exp\left(-\t\dt\sum_{k=1}^{\kappa
        }(\xi^i_{n+1,k})^4\right)(\xi^i_{n+1,\kappa})^p}{\sum_{j=1}^N\exp\left(-\t\dt\sum_{k=1}^{\kappa }(\xi^j_{n+1,k})^4\right)}\right)\leq \frac{e^{C_p\Dt}}{N}\sum_{i=1}^N\E((\xi^i_{n+1,0})^p)+(e^{C_p\dt}-1)\sum_{k=0}^{\kappa-1}e^{C_pk\dt}.$$
Inserting this bound in \eqref{decompmom}, one concludes that
\begin{align*}
 \E\left(\frac{\eta_{n+1}^N(gh)}{\eta_{n+1}^N(g)}\right) &\leq
 e^{C_p\Dt}\|h\|_{\kappa,p}\left(1+\E\left(\frac{1}{N}\sum_{i=1}^N(\xi^{i}_{n+1,0})^p\right)\right).
\end{align*}
For $n=0$, one deduces that $\E\left(\frac{\eta_{1}^N(gh)}{\eta_{1}^N(g)}\right)\leq
 e^{C_p\Dt}\|h\|_{\kappa,p}(1+\E(X_0^p))$, where $X_0$ is distributed according to the measure
$2\psi_I^2(x)1_{\{x>0\}}dx$.

For $n\geq 1$, since by a reasoning similar to the one made to obtain \eqref{condit},
$\ds{\E\left(\frac{1}{N}\sum_{i=1}^N(\xi^{i}_{n+1,0})^p\right)=\E\left(\frac{\eta_{n}^N(g(y)
      y_\kappa^p)}{\eta_n^N(g(y))}\right)}$,
 one also deduces that
$$\E\left(\frac{\eta_{n+1}^N(gh)}{\eta_{n+1}^N(g)}\right)\leq
e^{C_p\Dt}\|h\|_{\kappa,p}\E\left(\frac{\eta_{n}^N(g(1+y_\kappa^p))}{\eta^N_n(g)}\right).$$
The proof is completed by an obvious inductive reasoning.
\end{adem}

\begin{alem}\label{contmom}
For any $p\geq 2$, there is a constant $C_p$ such that 
$$\forall x\in\R_+^*,\;\forall t\geq 0,\;\E((X^x_t)^p)\leq
(1+x^p)e^{C_pt}-1,$$ 
where $X_t^x$ is defined by~(\ref{eq:sdex}). Therefore, if
$h:(\R_+^*)^\kappa\rightarrow \R$ is such that
$\|h\|_{\kappa,p}<+\infty$ then $\|Ph\|_{\kappa,p}\leq
e^{C_p\Dt}\|h\|_{\kappa,p}$, where the operator $P$ is defined by~(\ref{eq:P}).
\end{alem}
\begin{adem}
By Itô's formula, $d(X^x_t)^p=\left(\frac{p(p+1)}{2}(X^x_t)^{p-2}-\omega
  p(X^x_t)^{p}\right)dt+p(X^x_t)^{p-1}dW_t$. Hence
$$(X^x_t)^p\leq
x^p+\int_0^t\left(\frac{p(p+1)}{2}+\frac{p(p+1-2\omega)}{2}(X^x_s)^{p}\right)ds+p\int_0^t(X^x_s)^{p-1}dW_s.$$
Formally, taking expectations in this inequality, one obtains
\begin{align*}
   \E((X^x_t)^p)\leq x^p+\int_0^t\frac{p(p+1)}{2}+\frac{p(p+1-2\omega)}{2}\E((X^x_s)^{p})ds,
\end{align*}
and check by Gronwall's lemma that the conclusion holds with
$C_p=\frac{p(p+1)}{2}$. This formal argument can be made rigorous by a
standard localization procedure.\\
For $h:\R_+^\kappa\rightarrow \R$ such that $\|h\|_{\kappa,p}<+\infty$
one deduces that
$$\forall y\in\R_+^\kappa,\;|Ph(y)|\leq
\E|h(X^{y_\kappa}_\dt,\hdots,X^{y_\kappa}_{\kappa\dt})|\leq
C\|h\|_{\kappa,p}(1+\E((X^{y_\kappa}_{\kappa\dt})^p))\leq e^{C_p\Dt}\|h\|_{\kappa,p}(1+y_\kappa^p).$$\end{adem}
\begin{alem}\label{decroiss}
$$\forall(z_1,\hdots,z_N),(a_1,\hdots,a_N)\in\R_+^N\mbox{ with }\sum_{i=1}^Na_i>0,\;\forall p\geq 0,\;\forall c\geq
0,\;\frac{\sum_{i=1}^Na_iz_i^{p}e^{-cz_i^4}}{\sum_{i=1}^Na_ie^{-cz_i^4}}\leq
\frac{\sum_{i=1}^Na_iz_i^{p}}{\sum_{i=1}^Na_i}.$$
\end{alem}
\begin{adem}Let us set
   $f(c)=\frac{\sum_{i=1}^Na_iz_i^{p}e^{-cz_i^4}}{\sum_{i=1}^Na_ie^{-cz_i^4}}$. 
By Hölder's inequality, the derivative
$$f'(c)=\left(
  \frac{\sum_{i=1}^Na_iz_i^{p}e^{-cz_i^4}}{\sum_{i=1}^Na_ie^{-cz_i^4}}
  \frac{\sum_{i=1}^Na_iz_i^{4}e^{-cz_i^4}}{\sum_{i=1}^Na_ie^{-cz_i^4}} \right)-\frac{\sum_{i=1}^Na_iz_i^{p+4}e^{-cz_i^4}}{\sum_{i=1}^Na_ie^{-cz_i^4}}$$ 
is non positive. Hence for any 
$c\geq 0$, $f(c)\leq f(0)=\frac{\sum_{i=1}^Na_iz_i^{p}}{\sum_{i=1}^Na_i}$.
\end{adem}
\begin{alem}\label{contetag}
The sequence $(\eta_n(g))_{1\leq n\leq\nu}$ is bounded from below by a
positive constant non depending on~$\kappa$.\end{alem}
\begin{adem}
Since
$$\eta_n(g)=\frac{\E\left(\exp\left(-\t\dt\sum_{k=1}^{n\kappa}X_{k\dt}^4\right)\right)}{\E\left(\exp\left(-\t\dt\sum_{k=1}^{(n-1)\kappa}X_{k\dt}^4\right)\right)}\leq
1$$
the sequence $(\eta_n(g))_{1\leq n\leq\nu}$ is bounded from below by 
$$\prod_{n=1}^\nu
\eta_n(g)=\E\left(\exp\left(-\t\dt\sum_{k=1}^{\nu\kappa}X_{k\dt}^4\right)\right).$$
According to Lemma \ref{vitfort}, this expectation converges to
$\E\left(\exp\left(-\t\int_0^TX_{s}^4ds\right)\right)>0$ when $\kappa$
tends to $+\infty$, which concludes the proof.\end{adem}

We can now prove, as a corollary of Theorem~\ref{th:est_err}, the convergence of the
approximation $\overline{E_{\rm DMC}^{N,\nu,\kappa}}(T)$ of $E_{\rm
  DMC}(T)$, defined by:
\begin{equation}
\overline{E_{\rm DMC}^{N,\nu,\kappa}}(T)=\frac{3}{2}\omega
+\frac{\theta}{N}\sum_{i=1}^N (X^i_{\nu+1,0})^4.\label{expr3}
\end{equation}
\begin{acor}\label{cor:est_err}
$$
   \E\left|E_{\rm DMC}(T)-\overline{E_{\rm DMC}^{N,\nu,\kappa}}(T)\right|\leq
   \frac{C}{\nu\kappa}+\frac{C_\nu}{\sqrt{N}},
$$
where the constant $C$ only depends on $T$ and the constant $C_\nu$ on $T$ and $\nu$.
\end{acor}
\begin{adem}
By using the result of Theorem~\ref{th:est_err} and Cauchy Schwartz inequality, it is sufficient to prove
the estimate $\ds{\E\left(E_{\rm DMC}^{N,\nu,\kappa}(T)-\overline{E_{\rm
      DMC}^{N,\nu,\kappa}}(T)\right)^2\leq \frac{C_\nu}{N}}$. Let us
denote in this proof $\xi^i_{\nu+1,0}=X^i_{\nu+1,0}$ for $1\leq i \leq N$. We have:
\begin{equation*}
E_{\rm DMC}^{N,\nu,\kappa}(T)-\overline{E_{\rm DMC}^{N,\nu,\kappa}}(T)
= \theta\left(\frac{\eta_\nu^N(g \, y_\kappa^4)}{\eta_\nu^N(g)} -
  \frac{1}{N} \sum_{i=1}^N (\xi^i_{\nu+1,0})^4\right)
= \theta\left( \E\left( \frac{1}{N} \sum_{i=1}^N (\xi^i_{\nu+1,0})^4 \bigg|
    \xi_\nu \right) -
  \frac{1}{N} \sum_{i=1}^N (\xi^i_{\nu+1,0})^4\right)
\end{equation*}
by using the fact that, for any function $f:\R_+^* \to \R_+$, 
\begin{equation}\label{eq:rec}
\E\left( \frac{1}{N}
  \sum_{i=1}^N f(\xi^i_{\nu+1,0}) \bigg| \xi_\nu \right)=\frac{\eta_\nu^N(g(y)
  \,f(y_\kappa))}{\eta_\nu^N(g(y))},
\end{equation} which is obtained by a reasoning similar to the
one made to prove~\eqref{condit}. Now, using the same method as to
obtain~\eqref{majot}, one easily gets the estimate:
$$\E\left( E_{\rm DMC}^{N,\nu,\kappa}(T)-\overline{E_{\rm
      DMC}^{N,\nu,\kappa}}(T) \right)^2 \leq \frac{\theta^2}{N}
\E\left(\frac{1}{N}\sum_{i=1}^N (\xi^i_{\nu+1,0})^8\right)= \frac{\theta^2}{N} \E\left(
\frac{\eta_\nu^N(g(y)
  \,(y_\kappa)^8)}{\eta_\nu^N(g(y))}\right),
$$
by using again~(\ref{eq:rec}). Lemma~\ref{majomom} completes the proof.
\end{adem}

We end this Section by proving that Proposition~\ref{propvitpart} also
holds for the numerical scheme~\eqref{schemaurel}.
\begin{aprop}\label{prop:Markov}
Let us consider the Markov chain $(\bar{X}_{j\delta t})_{0\leq j\leq K}$
generated by the explicit scheme \eqref{schemaurel} and denote by $Q$ its
transition kernel. We now define the measure
$\eta_n$ by replacing $(X_{j\delta t})_{0\leq j\leq
  K}$ with $(\bar{X}_{j\delta t})_{0\leq j\leq K}$ in~(\ref{eq:eta_n}), and we define accordingly the
   evolution of the particle system: conditionally on $\xi_n$,
 the vectors $(X^i_{n+1,0},X^i_{n+1,\dt},\hdots,X^i_{n+1,\kappa\dt})_{1\leq i\leq N}$
   are independent, with $(X^i_{n+1,0})_{1 \le i \le N}$ distributed according to
 the selection algorithm (S1) (see~\eqref{messel}), (S2)
 (see~\eqref{messel2}) or (S3) (see~\eqref{eq:strat}), and $(X^i_{n+1,j\dt})_{0\leq j\leq \kappa}$ a Markov
   chain with transition kernel $Q$. Then, we have:
$$
   \E\left|\frac{\eta^N_\nu(gy_\kappa^4)}{\eta^N_\nu(g)}-\frac{\eta_\nu(gy_{\kappa}^4)}{\eta_\nu(g)}\right|\leq \frac{C_\nu}{\sqrt{N}}.
$$
\end{aprop} 
\begin{adem}
Looking carefully at the proof
   of Proposition \ref{propvitpart} above, one remarks that \eqref{vitpart}
   holds in this framework as soon as
   Lemma~\ref{contetag} holds, and the following property, which
   replaces Lemma \ref{contmom}, is satisfied:
\begin{equation}\label{eq:hyp_schemaurel}
\exists C>0,\;\forall x\in\R_+,\;Qf(x)\leq e^{C\delta
  t}(1+f(x))-1\mbox{ for }f(x)\equiv x^4\mbox{ and }f(x)\equiv x^8.
\end{equation}
Let us first prove~(\ref{eq:hyp_schemaurel}). We have:
$Qf(x)=\E\left(f\left(\bar{X}^x_{\delta t}\right)\right)$ where $\bar{X}^x_\dt=\left((1-\omega
    \delta t)^2x^2+2xW_\dt+\frac{W_\dt^2}{(1-\omega
    \delta t)^2}+2\delta t\right)^{1/2}$. Now, for $q\in\N^*$,
$$(\bar{X}^x_{\delta t})^{2q}=\sum_{j_1+j_2+j_3=q}\frac{q!}{j_1!j_2!j_3!}\,(1-\omega\dt)^{2j_1}\,2^{j_2}\,x^{2j_1+j_2}\,W_\dt^{j_2}\left(\frac{W_\dt^2}{(1-\omega
    \delta t)^2}+2\delta t\right)^{j_3},$$
where the indices $(j_1,j_2,j_3)$ are non negative integers. Remarking that the expectation of the terms with $j_2$ odd vanishes and
    then using Young's inequality, one deduces that for $\dt\leq
    \frac{1}{2\omega}$,
\begin{eqnarray}
   \E\left((\bar{X}^x_{\delta t})^{2q}\right)&\leq& (1-\omega\dt)^{2q}x^{2q}+\E\left(\left(\frac{W_\dt^2}{(1-\omega
    \delta t)^2}+2\delta
    t\right)^{q}\right)+C_q \!\!\!\!\!\!\!\!\!\sum_{\stackrel{j_1+j_2+j_3=q}{j_1<q,j_2\;even\;,j_3<q}}\!\!\!\!\!\!\!\!\!x^{2\left(q-\frac{j_2+2j_3}{2}\right)}\dt^{\frac{j_2+2j_3}{2}},\nonumber\\
&\leq& x^{2q}+C_q\dt+C_q\!\!\!\!\!\!\!\!\!\sum_{\stackrel{j_1+j_2+j_3=q}{j_1<q,j_2
    \;even\;,
    j_3<q}}\!\!\!\!\!\!\!\!\!\left(x^{2q}\dt+\dt^{1+q\left(1-\frac{2}{j_2+2j_3}\right)}\right),\nonumber\\
&\leq& (1+C_q\dt)x^{2q}+C_q\dt\leq e^{C_q\dt}(1+x^{2q})-1. \label{eq:Xbar_estim}
\end{eqnarray}
Let us now prove Lemma~\ref{contetag} for the scheme
\eqref{schemaurel}. As noticed in the proof of  Lemma~\ref{contetag}
above, it is sufficient to bound from below
$\E\left(\exp\left(-\t\dt\sum_{k=1}^{\nu\kappa}\bar{X}_{k\dt}^4\right)\right)$.
By Jensen inequality, we have
$\E\left(\exp\left(-\t\dt\sum_{k=1}^{\nu\kappa}\bar{X}_{k\dt}^4\right)\right)
 \geq \exp\left(-\t\frac{T}{\nu \kappa}\sum_{k=1}^{\nu\kappa}
   \E\left(\bar{X}_{k\dt}^4\right)\right)$. By
 using~(\ref{eq:Xbar_estim}), it is easy to prove by induction that
 $\E\left(\bar{X}_{k\dt}^4\right)\leq e^{C_2 k \dt}
 (1+\E\left(\bar{X}_{0}^4\right))-1$ and this concludes the proof of
 Lemma~\ref{contetag} in this framework.
\end{adem}
In order to obtain a complete convergence result of the
form~(\ref{eq:est_err}) for the scheme \eqref{schemaurel}, it remains to
prove the complementary bound \eqref{discretaurel}, that we have not
obtained so far. However, we will check by numerical simulations that (\ref{eq:est_err}) still holds.

\section{Numerical results}\label{sec:num_res}

\subsection{Computation of a reference solution by a spectral method}\label{sec:num_res_spec}

In this section, we would like to explain how we can obtain a very
precise reference solution by using a partial differential equation approach to compute
$E_{\rm DMC}(T)$ (see~\cite{CJL}).

\subsubsection{A partial differential equation approach to compute $E_{\rm DMC}(T)$}

Let us introduce the solution $\phi$ to the following partial
differential equation for :
\begin{equation}\label{eq:phi}
\left\{
\begin{array}{l}
\displaystyle{\frac{\partial \phi}{\partial t}=-H \phi},\;(t,x) \in \R_+ \times \R\\
\phi(0,x)=\psi_I(x),\;x\in\R
\end{array}
\right.
\end{equation}
where $H$ (resp. $\psi_I$) is defined by~(\ref{eq:H})
(resp.~(\ref{eq:psi_I})). Since $\psi_I \in {\mathcal H}$, it
is a standard result that this problem admits a unique solution
$\phi  \in C^0(\R_+,{\mathcal H}) \cap C^0(\R_+^*,D_{\mathcal H}(H))
\cap C^1(\R_+^*,{\mathcal H})$. The function $\phi$ is
regular and odd, and therefore is such that $\phi(t,0)=0$ for all $t
\geq 0$. Therefore the function $\phi$ is also solution to the following
partial differential equation:
\begin{equation}\label{eq:phi_prime}
\left\{
\begin{array}{l}
\displaystyle{\frac{\partial \phi}{\partial t}=-H \phi},\;(t,x) \in \R_+ \times \R\\
\phi(t,0)=0,\; t\geq 0\\
\phi(0,x)=\psi_I(x),\;x\in\R.
\end{array}
\right.
\end{equation}
In~\cite{CJL}, we have shown that since $\phi$ satisfies~(\ref{eq:phi_prime}), we can express $E_{\rm DMC}(t)$
(defined by~(\ref{eq:EDMC}))
using the function $\phi$ (see Proposition 11 in~\cite{CJL}):
\begin{equation}\label{eq:EDMC_phi}
E_{\rm DMC}(t)=\frac{\langle H \psi_I,\phi(t)\rangle}{\langle
  \psi_I,\phi(t)\rangle}.
\end{equation}
Our reference solution $E_{\rm DMC}(T)$ will rely on formula
\eqref{eq:EDMC_phi} after discretization of (\ref{eq:phi}) by a spectral
method.


\subsubsection{Computation of the wave function $\phi$}

We will briefly present the spectral method developed to compute an
approximation of $\phi$.
We recall that the Hermite polynomials are defined by  :
$$
\forall n \in \N,\; h_n(x)=(-1)^n e^{x^2}\displaystyle\frac{d^n}{dx^n}(e^{-x^2}).
$$ 
We introduce the eigenfunctions 
of the operator $H_0$, 
normalized for the $L^2(\R)$ norm associated with the eigenvalues
$E_n=\omega(n+1/2)$ for $n \geq 0$,
$$
\varphi_n(x)=\displaystyle h_n(\sqrt \omega x) \exp({-\frac{1}{2}\omega x^2})
\left(\frac{(\omega /\pi)^{1/4}}{\sqrt{2^n n!}}\right).
$$
It is well known
that the vector space spanned by the set of functions
$\{\varphi_{2k+1}\}_{k \geq 0}$ is dense in ${\mathcal V}_0=\{\varphi\in
H^1(\R) \cap {\mathcal H} \;| \;  x\varphi \in L^2\}$, which is the
domain of the quadratic form associated with $H_0$.

Let us now introduce the functional space  
 ${\mathcal V}=\{\varphi\in H^1(\R) \cap {\mathcal H} \;| \;
 x^2\varphi\in L^2\}$, which is the
domain of the quadratic form associated with $H$. The set of functions $\{\varphi_{2k+1}\}_{k \geq 0}$ is also a
 basis of ${\mathcal V}$.

 Let ${\mathcal V}_n
 =Span({\varphi_1},{\varphi_3},\dots ,{\varphi}_{2n-1})$. We use this
 approximation space to build the following Galerkin scheme
 for~(\ref{eq:phi}): find $\phi_n \in C^0(\R_+,{\mathcal V}_{n})$ such
 that\footnote{Notice that $\psi_I=\varphi_1\in {\mathcal V}_n$.}  
$\phi_n(0,x)=\psi_I$, and $\forall \varphi \in {\mathcal V}_n$
\begin{equation}\label{eq:1}
\displaystyle\left\langle\frac {\partial \phi_n(t)}{\partial t},\varphi\right\rangle=
\displaystyle-\left\langle H \phi_n(x,t),\varphi\right\rangle.
\end{equation}

We diagonalize the operator $H$ restricted to ${\mathcal V}_n$. We denote 
$({\varphi_0^n},{\varphi}_2^n,\dots ,{\varphi}_{n-1}^n)$ 
the eigenfunctions and ${E_0^n},{E_2^n},\dots, { E_{n-1}^n}$
 the associated eigenvalues. Because of the symmetry of $H$, it is easy
 to check that ${\mathcal V}_n$ can also be spanned by $({\varphi_0^n},{\varphi}_2^n,\dots ,{\varphi}_{n-1}^n)$:
\begin{equation}\label{eq:span}
{\mathcal V}_n=Span({\varphi_0^n},{\varphi}_2^n,\dots ,{\varphi}_{n-1}^n).
\end{equation}
Since for $t\geq 0$, $\phi_n(t,.) \in {\mathcal V}_n$, there exists $u_{k}(t), \;k=0,\dots,n-1$, such that 
\begin{equation}\label{eq:psi}
\displaystyle\phi_n=\sum_{k=0}^{n-1} u_{k}(t){\varphi_{k}^n}.
\end{equation}
In view of (\ref{eq:span}) and (\ref{eq:psi}), (\ref{eq:1}) is equivalent
 to the  equations: $\forall i=0,\dots,n-1$,
$$
\begin{array}{ccl}
\displaystyle\sum_{k=0}^{n-1}\frac {\partial u_{k}(t)}{\partial t} 
\left\langle{\varphi}_{k}^n,{\varphi}_{i}^n\right\rangle&=&
\displaystyle-\left\langle H\sum_{k=0}^{n-1} u_{k}(t) 
{\varphi}_{k}^n,{\varphi}_{i}^n\right\rangle, \\
&=&\displaystyle-\sum_{k=0}^{n-1}{ E}_{k}^nu_{k}(t)\left\langle{\varphi}_{k}^n,  
{\varphi}_{i}^n\right\rangle.
\end{array}
$$
We deduce that  $\forall k =0,\dots,n-1,$
$$
\displaystyle\frac {\partial u_{k}(t)}{\partial t} = -{E}_{k}^n u_{k}(t),
$$
so that
\begin{equation}\label{eq:psi2}
\displaystyle\phi_n(t,x)=\sum_{k=0}^{n-1}u_{k}(0)\exp(-{E}_{k}^n t){\varphi_{k}^n}(x),
\end{equation}
where $u_{k}(0)=\left\langle\psi_I,{\varphi}_k^n\right\rangle$.

\begin{arem}
The eigenfunctions of $H$ are obtained by diagonalization of the matrix
$A=(a_{ij})_{i,j=0,\dots,n-1}$ with $\:\:\forall i,j=0,\dots,n-1$ :
$$
\begin{array}{ccl}
a_{ij}&=&\left\langle H\varphi_{2i+1},\varphi_{2j+1}\right\rangle, \\
           &=&\left\langle H_0\varphi_{2i+1},\varphi_{2j+1}\right\rangle+
             \theta \left\langle x^4\varphi_{2i+1},\varphi_{2j+1}\right\rangle,\\
  &=&\delta_{ij} \, \omega \, (2i+\frac{3}{2}) + \theta\, \left\langle  x^4\varphi_{2i+1},
\varphi_{2j+1}\right\rangle.
\end{array}
$$
We can use the  n--point Gauss-Hermite formula to  deal with the integration of the
 second term on the right. We recall that this method  provides an exact
 result for
$\int_{-\infty}^{+\infty}p(x) \exp(-x^2)dx $ as long as $p$ is a
polynomial of degree $2n-1$ or less.
 
\end{arem}
\subsubsection{Approximation of  $E_{\rm DMC}(T)$}

We now use formula~(\ref{eq:EDMC_phi}) to approximate $E_{\rm DMC}(T)$.
By an elementary calculation, we obtain the following approximation:
\begin{equation}\label{eq:E2}
E_{\rm DMC}(T) \simeq \displaystyle\frac{{E}_0^n+\displaystyle\sum_{i=1}^{n-1}\displaystyle\frac{u_{i}(0)\left\langle{\varphi}_i^n,\varphi_1\right\rangle}
{u_{0}(0)\left\langle{\varphi}_1^n,\varphi_1\right\rangle}{E}_{i}^n\exp(-({E}_{i}^n-{E}_0^n)T)}
{1+\displaystyle\sum_{i=1}^{n-1}\frac{u_{i}(0)\left\langle{\varphi}_i^n,\varphi_1\right\rangle}
{u_{0}(0)\left\langle{\varphi}_1^n,\varphi_1\right\rangle}\exp(-({E}_{i}^n-{E}_0^n)T)}.
\end{equation}
In our test cases, we have observed that $n=40$ is enough to reach
convergence.

Notice that for a given $n$, the convergence in time to the lowest eigenvalue
${E}_0^n$ is exponentially fast, with an exponent equal to the 
spectral gap ${E}_1^n-{E}_0^n$.

\subsection{Numerical results of Monte Carlo simulations}

In this section, we perform various numerical experiments to validate
our theoretical results, and to explore some features of DMC
computation. In particular, we propose in Section~\ref{sec:optim_var} an
empirical method to determine the optimal number of reconfigurations. In
all the computations, the final time is $T=5$, which appears to be
sufficiently large for the convergence $t \to \infty$ to be achieved
with enough accuracy.

\subsubsection{Error and variance as a function of the numerical parameters}\label{sec:err_var}

\begin{figure}
\begin{center}
\psfrag {e(dt)}{$e(\delta t)$}
\psfrag{v(dt)}{$v(\delta t)$}
\psfrag{e(N)}{$e(N)$}
\psfrag{v(N)}{$v(N)$}
\psfrag{e(Nu)}{$e(\nu)$}
\psfrag{v(Nu)}{$v(\nu)$}
\begin{tabular}{cc}
      \includegraphics[width=6cm]{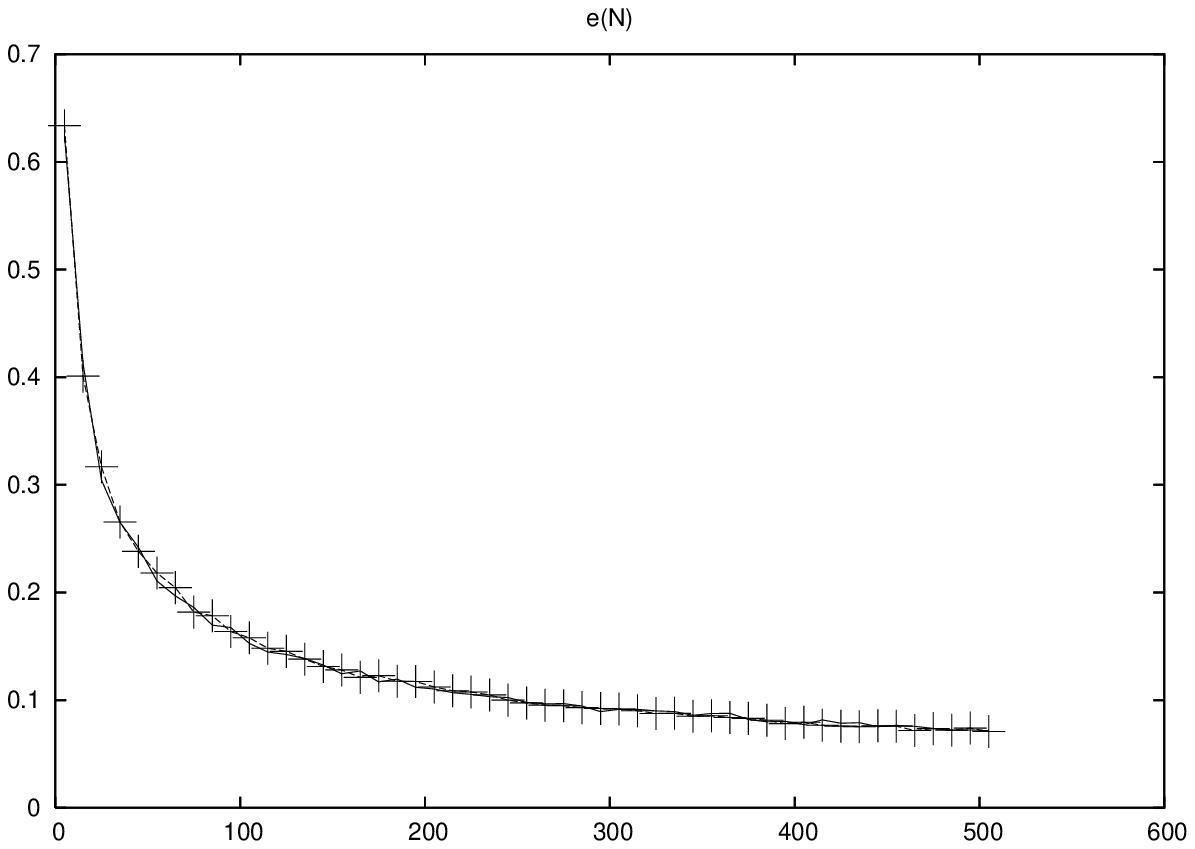} &
      \includegraphics[width=6cm]{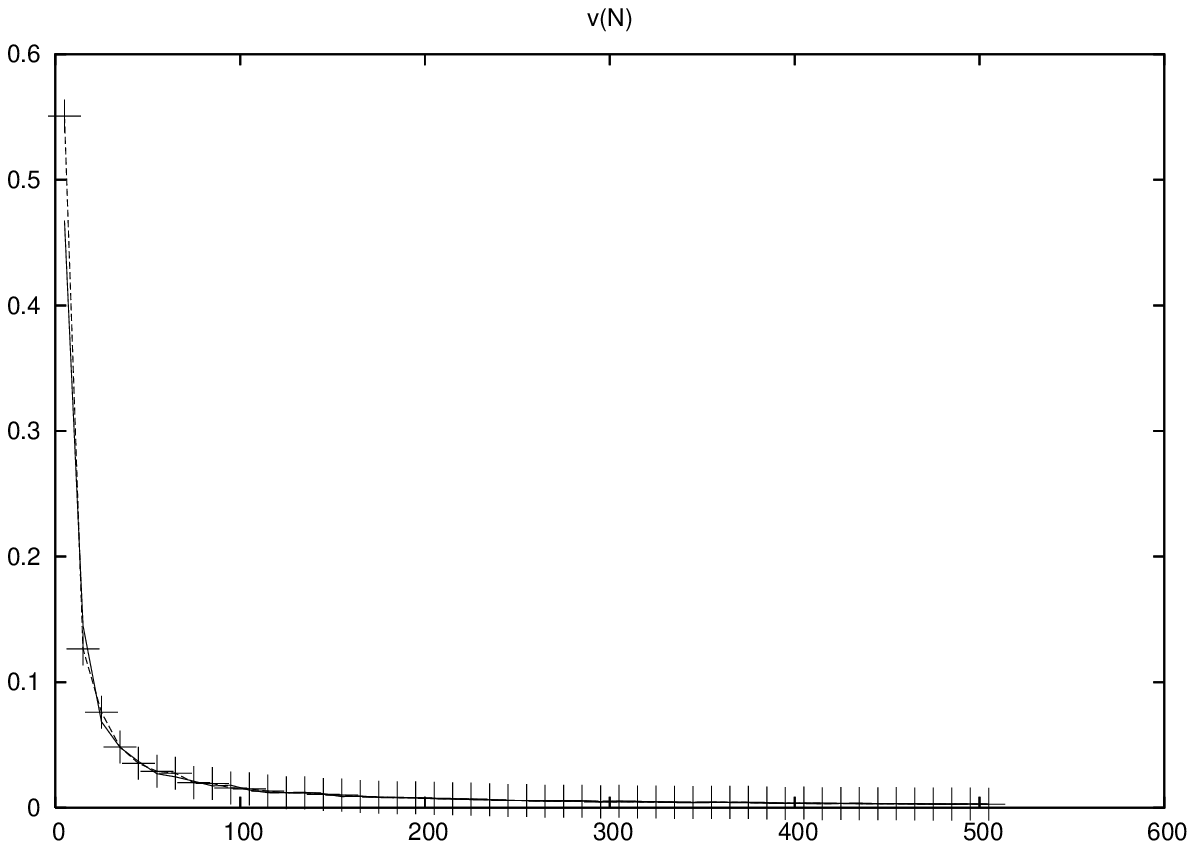} \\
      \includegraphics[width=6cm]{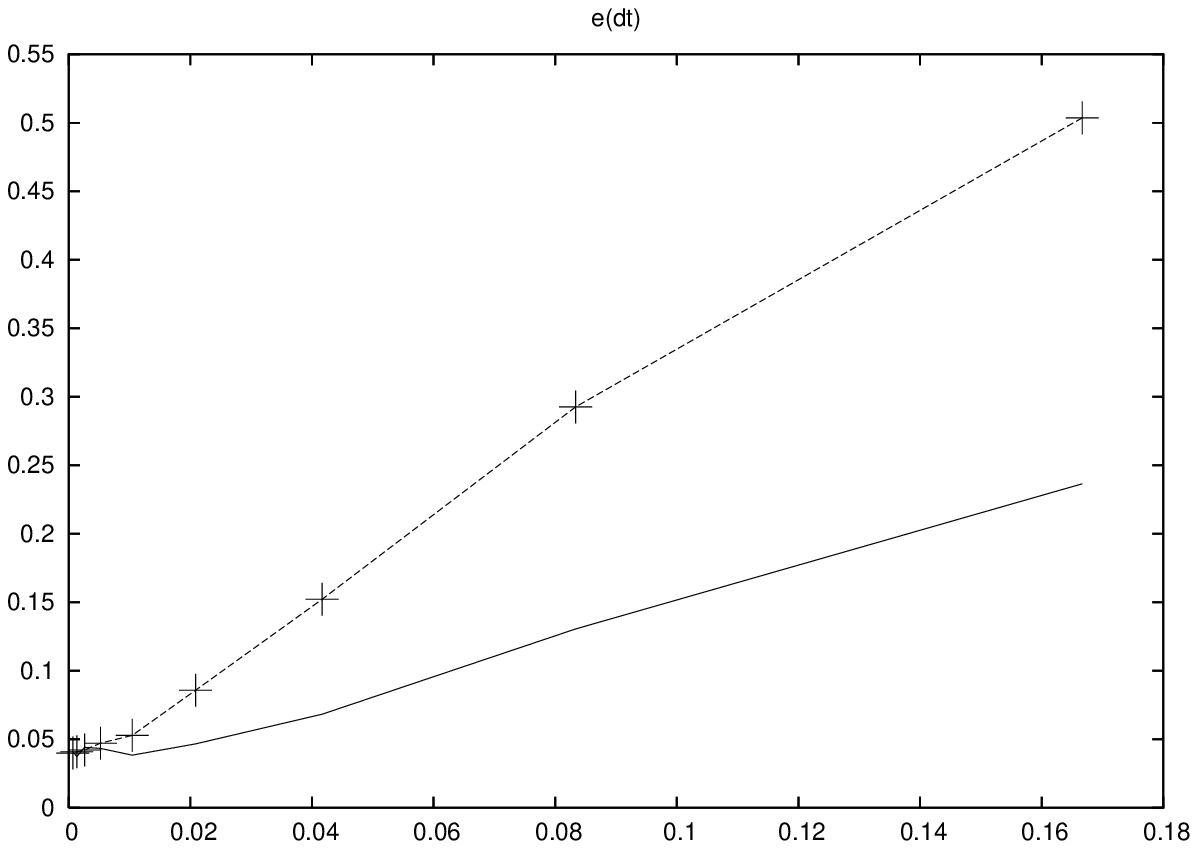} &
      \includegraphics[width=6cm]{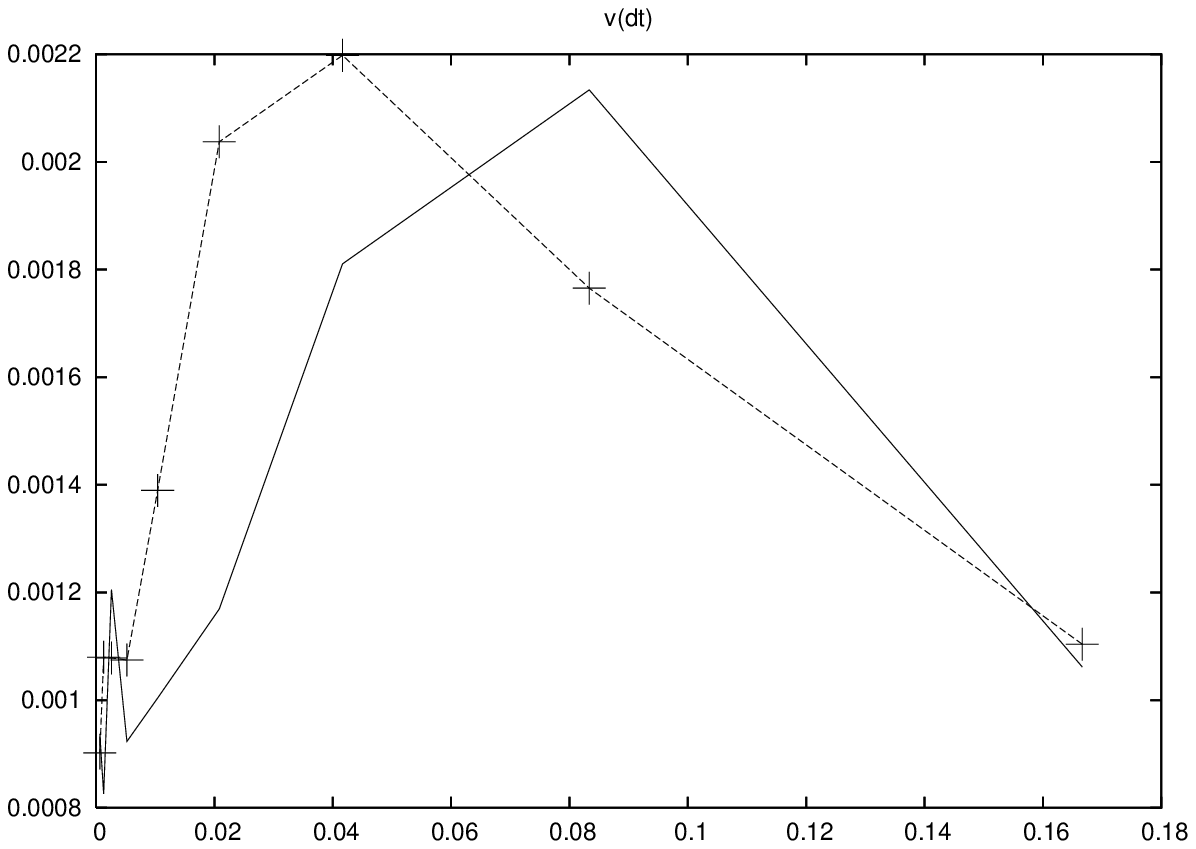} \\
      \includegraphics[width=6cm]{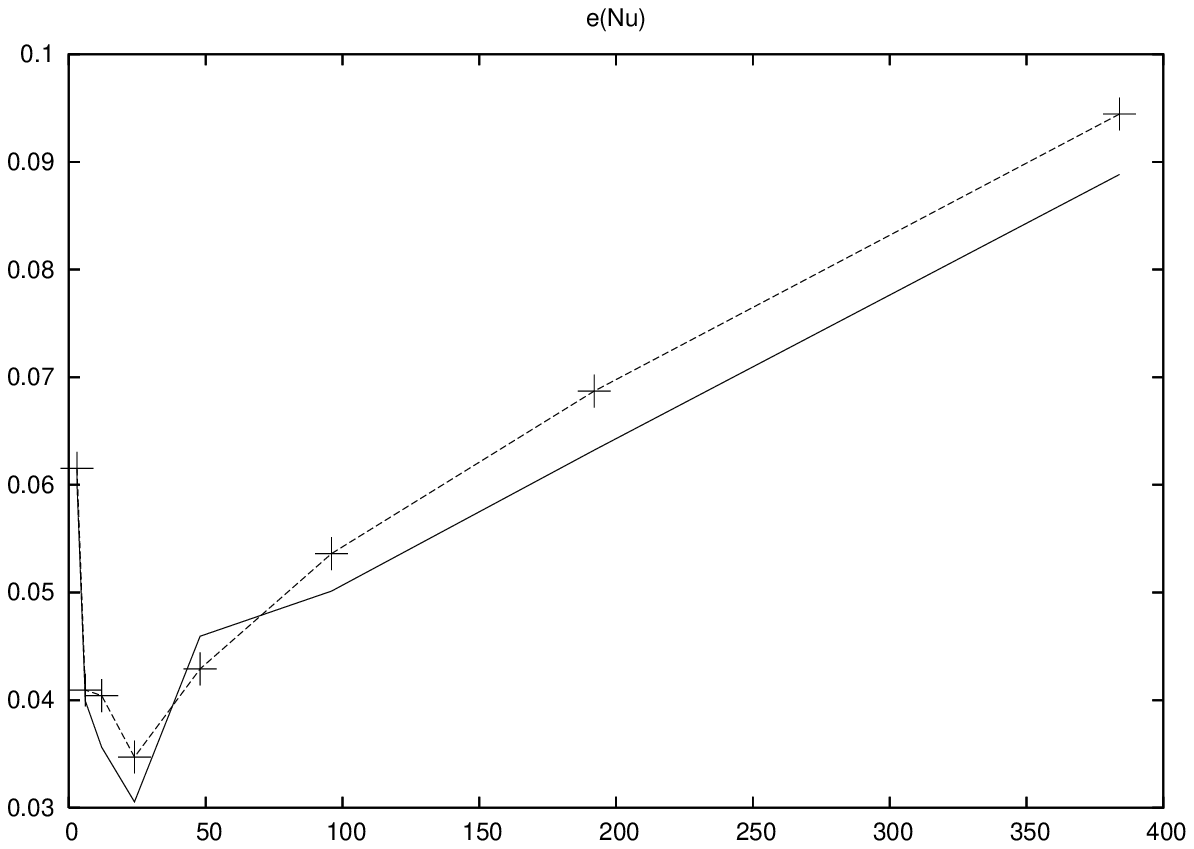} &
      \includegraphics[width=6cm]{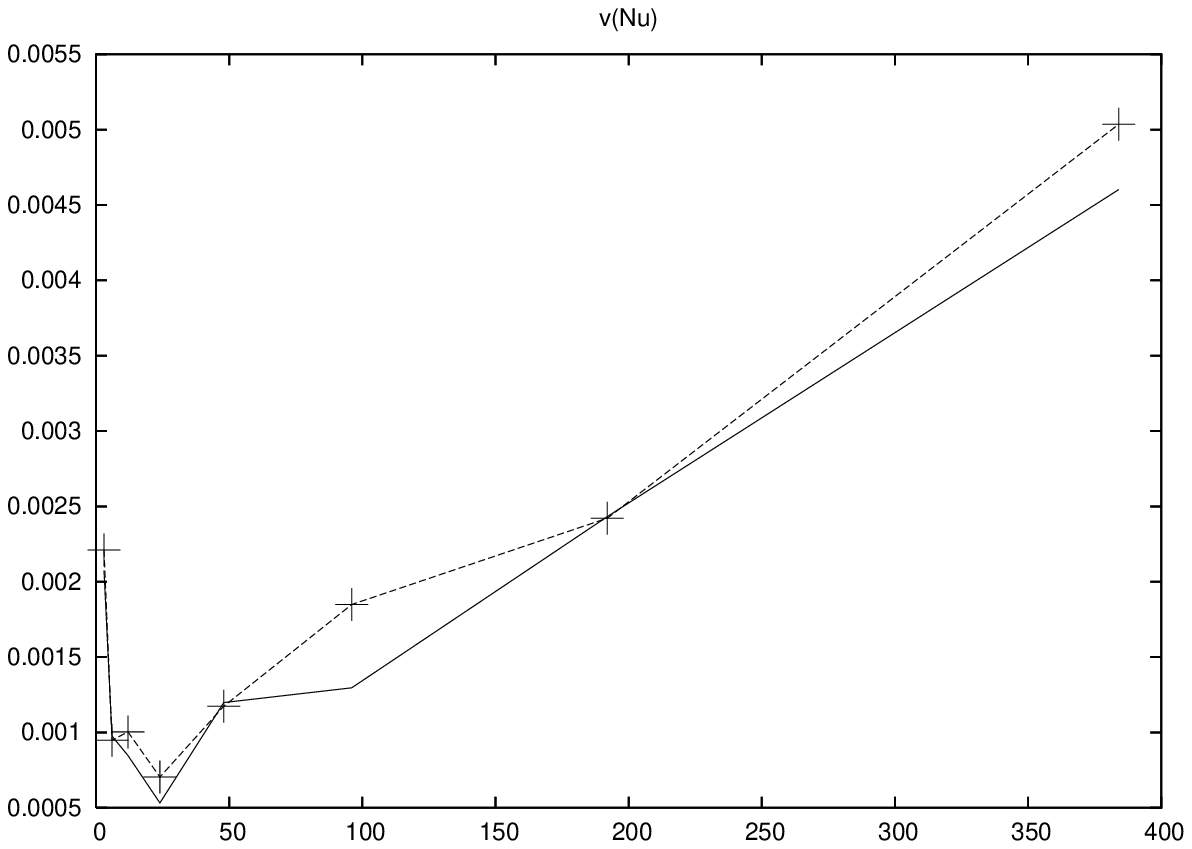} \\
    \end{tabular}
\caption{Expectation and variance of the error when \eqref{eds} is
  discretized according to the method 
  described in Appendix (dotted curve) and according to the
  scheme~(\ref{schemaurel}) (solid curve).}\label{fig}
  \end{center}
\end{figure}


We represent on Figure \ref{fig}, the expectation $e$ and the variance
$v$ of the error : $\left|E_{\rm DMC}^{N,\nu,T / (\nu \delta t)}(T)-E_{\rm DMC}(T) \right|$ as a function of the  
number of walkers $N$, the time step $\delta t$ and the number of
reconfigurations $\nu-1$, where $E_{\rm DMC}(T)$ is approximated
using~(\ref{eq:E2}) and $E_{\rm DMC}^{N,\nu,T / (\nu \delta t)}(T)$ is defined
by~(\ref{expr2}). The multinomial resampling method (which is (S1) with
  $\epsilon_n=0$) was used.

The top figures represent the expectation of the error and its  variance 
according to the number of walkers. 
To compute these quantities, we perform $2000$ independent realizations,
 with the  number of reconfigurations 
$\nu -1=50$, a small time step $\delta t=5.10^{-3}$ and  $\theta=0.5$.
The simulations confirm the theoretical result : the
error decreases as $C/\sqrt{N}$.

The effect of the time step is shown on the two figures in the center.
The numerical parameters are: a large number of particles $N=5000$, number of configurations
 $\nu -1=30$, $\theta=2$ and $300$ independent realizations. We can  see on the figure on the left that the  error 
decreases linearly as the time step decreases. We also remark that the
error is smaller with the approximate scheme~(\ref{schemaurel}) than when using
the exact simulation of the SDE~(\ref{eds}) proposed in the
Appendix. This rather amazing result can be interpreted as follows. When
using the exact simulation of the SDE, there is only one source of error
related to the time discretization, namely the approximation of the
integral in the exponential factor in~(\ref{eq:EDMC}). When using the
scheme~(\ref{schemaurel}), we add a weak error term which seems to
partly compensate the previous one.

The last figures represent the effect of the number of
  reconfiguration steps. The numerical parameters are: time step $\delta
  t=5.10^{-3}$, number of particles $N=5000$, $\theta=2$ and $300$ independent realizations.
The curve representing the variation of the error according to the number of
 reconfigurations has the shape of a basin.
We deduce that on the one hand a small number of reconfigurations has 
the disadvantage that walkers with increasingly differing weights are kept. On the other hand 
a large number of reconfigurations introduces much noise. An optimal
number of reconfiguration seems to lie between 20 and 50.

\subsubsection{Optimal number of reconfigurations}\label{sec:optim_var}

\begin{figure}[ht]
\begin{center}
\psfrag{var(EDMC)}{$\tilde{v}(\nu)$}
\includegraphics[height=6cm]{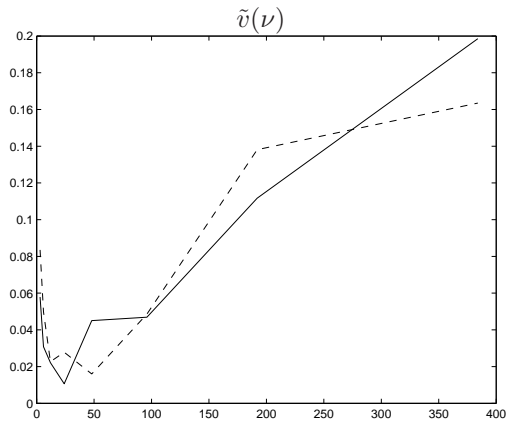}
\caption{Variance of $E_{\rm DMC}^{N,\nu,T / (\nu \delta t)}(T)$ in function of the number of
  reconfigurations when \eqref{eds} is
  discretized according to the method 
  described in Appendix (solid curve) and according to the
  scheme~(\ref{schemaurel}) (dashed curve).}\label{fig:var_EDMC}
\end{center}
\end{figure}

On Figure~\ref{fig:var_EDMC}, we check that the
optimal number of reconfigurations in terms of the variance $\tilde{v}$
of $E_{\rm DMC}^{N,\nu,T / (\nu \delta t)}(T)$ (and not of the error as
in Section~\ref{sec:err_var}) is also obtained
for a number of reconfiguration which seems to lie between 20 and 50
(using again the multinomial resampling method). The numerical
parameters are those considered for the figures below in
Figure~\ref{fig}: time step $\delta
  t=5.10^{-3}$, number of particles $N=5000$, $\theta=2$ and $300$
  independent realizations. We have not studied how the optimal number of reconfigurations varies according to the other
numerical parameters.

\begin{figure}[ht]
\begin{center}
\psfrag{var EDMC(t)}{$\tilde{v}(t)$}
\includegraphics[height=6cm]{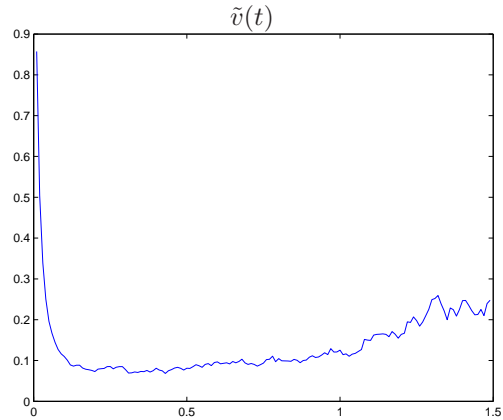}
\caption{Variance of $E_{\rm DMC}^{N,1,t / \delta t}(t)$ as a function of time $t$.}\label{Nreconf}
\end{center}
\end{figure}

We have investigated a practical method to estimate numerically the optimal number of 
reconfigurations. On Figure~\ref{Nreconf} we represent the variance of
$E_{\rm DMC}^{N,1,t / \delta t}(t)$ according 
to time $t$, without any reconfiguration step (which corresponds to $\nu=1$). The other numerical
parameters are again those considered for the figures below in
Figure~\ref{fig}. We observe that the variance is minimal at $t^*
\approx 0.25 $. We remark that $\nu =
T/t^*= 20$ is close to the optimal number of reconfigurations 
obtained on the previous figures. We have checked this empirical result
for various sets of the parameters. 
It seems that the optimal
number of reconfigurations is related to $T/t^*$ where $t^*$ minimizes
the variance of $E_{\rm DMC}^{N,1,t / \delta t}(t)$. Since $\nu=1$, no
selection step occurs and the particles are thus
independent. According to the multidimensional
central limit theorem, the variance of $E_{\rm DMC}^{N,1,t / \delta
  t}(t)$ can be approximated by
$$\frac{1}{N}\left( \frac{{\rm Var}(Y_t)}{ (\E(Z_t))^2} - 2 \E(Y_t) \frac{{\rm
       Covar}(Y_t,Z_t)}{(\E(Z_t))^3} + (\E(Y_t))^2 \frac{{\rm Var}(Z_t)}{(\E(Z_t))^4}\right)$$
where
$$Y_t=E_L(X_t) \exp\left(-\dt \sum_{k=1}^{t /\dt} E_L(X_{k \dt}) \right)$$
and
$$Z_t=\exp\left(-\dt \sum_{k=1}^{t /\dt} E_L(X_{k \dt}) \right).$$

 Therefore, the optimal number of reconfiguration steps
could be estimated by this method, through a precomputation over a
few independent trajectories.

 \subsubsection{Comparison of the resampling algorithms}

\begin{figure}[ht]
\begin{center}
\begin{tabular}{cc}
      \includegraphics[width=8cm]{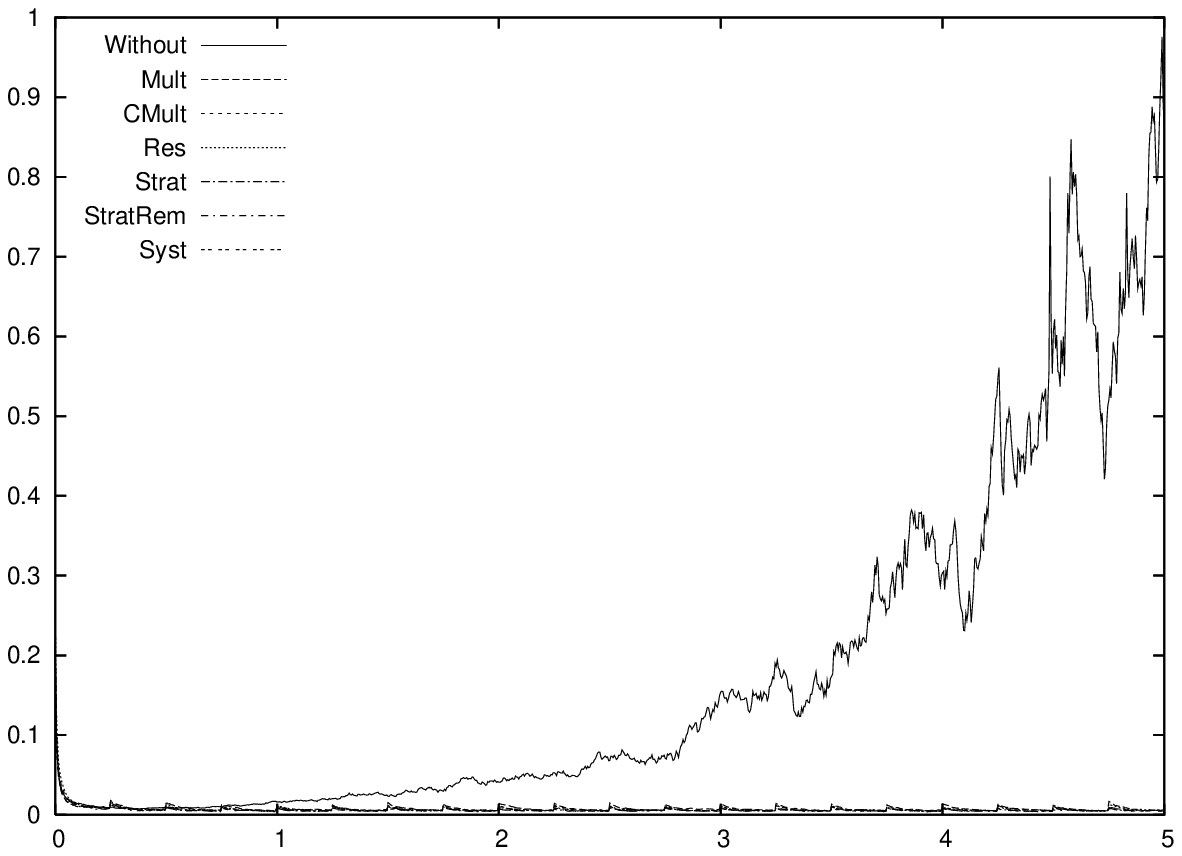} &
      \includegraphics[width=8cm]{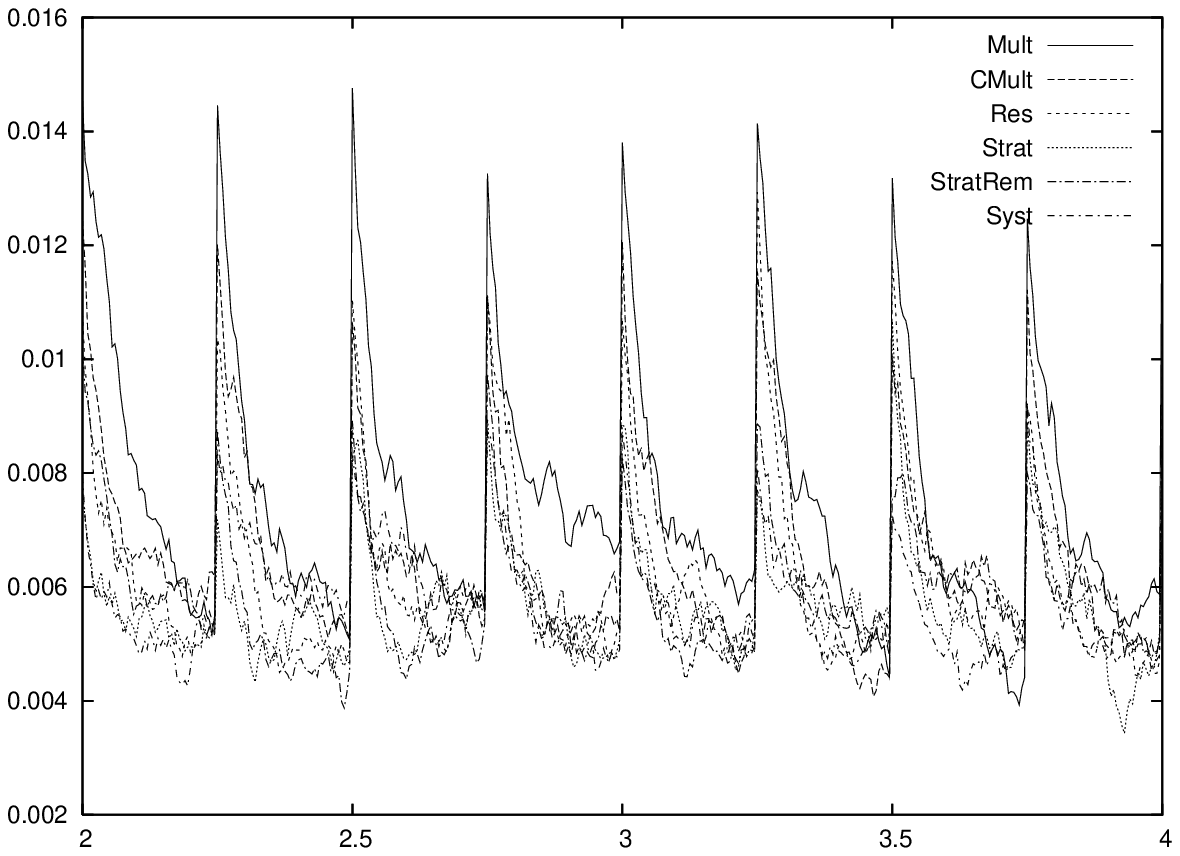}
\end{tabular}
\caption{Variance of $E_{\rm DMC}^{N,\nu,T / (\nu \delta t)}(t)$ as a
  function of time $t$, for various resampling algorithms: Without =
  without resampling, Mult = multinomial resampling, CMult = correlated
  multinomial resampling, Res = residual resampling, Strat = stratified
  resampling, StratRem = stratified remainder resampling, Syst = systematic resampling.}\label{fig:var_resamp}
  \end{center}
\end{figure}

We finally compare various resampling algorithms on
Figure~\ref{fig:var_resamp}, where the variance of $E_{\rm DMC}^{N,\nu,T / (\nu \delta t)}(t)$ as a
  function of time is represented. The numerical parameters are:
  $N=1000$, $\dt=5.10^{-3}$, $\nu - 1 =20$, $\theta=2$ and 200
  independent realizations.

We first observe on the figure on the
  left that without any resampling, the variance of the results explodes
  with increasing time. This shows the necessity to use resampling
  algorithms. We compare the following resampling algorithms: multinomial resampling (which is (S1) with
  $\epsilon_n=0$), correlated multinomial resampling (which is (S1) with
  $\epsilon_n=1 / \max_{1\leq i\leq N}g(\xi_n^i)$), residual resampling
  (which is (S2)), stratified resampling (which is (S3)),
  stratified remainder resampling (which combines residual and
  stratified resampling, see Remark~\ref{rem:SRR}) and
  systematic resampling (which corresponds to stratified resampling with
  $U^1_n= \ldots=U^N_n=U_n$, see the Introduction). We observe that, as
  expected, when more correlation is introduced, the variance due to the
  resampling is reduced. The multinomial resampling method is generally
  the worse, while the best resampling methods seem to be systematic
  resampling or stratified remainder resampling.

\section*{Conclusion}

In this paper, we have proved on a simple example convergence of
numerical implementations of the DMC method with a fixed number of
walkers. The theoretical rates of convergence are confirmed by
numerical experiments and are likely to hold in more general situations. We
have also checked numerically the existence of an optimal number of
reconfiguration steps. Various resampling algorithms have been
considered, both theoretically and numerically. For future work, we plan to investigate criteria
devoted to the choice of the number of reconfiguration steps. One
interesting direction is the use of automatic criteria based on a
measure of the discrepancy between the weights carried by the walkers to
decide when to perform a reconfiguration step.

\section*{Appendix : Simulation of the stochastic differential equation
  \eqref{eds}}
In this appendix, we show that it is possible to simulate exactly in law
the $(K+1)$-plet
$(X_0,X_{\dt},\hdots,X_{K \dt})$, where $X_t$ is defined by~(\ref{eds}). Let $(G,U)$ denote a couple of independent random variables with $G$
normal and
$U$ uniformly distributed on the interval $[0,1]$.
\subsection*{Simulation of the increment $X_t-X_s$, for $t \geq s$.}
The square $R_t$ of the norm of a $3$-dimensional Brownian motion $\uu{W}_t=\left(\uu{W}^1_t,\uu{W}^2_t,\uu{W}^3_t\right)$ solves
$dR_t=3dt+2\sqrt{R_t}dB_t$ where $\ds{B_t=\int_0^t \frac{\uu{W}_s \cdot
  d\uu{W}_s}{\|\uu{W}_s \|}}$ is a one-dimensional Brownian
motion. Hence $\ds{\rho_t=\frac{R_t}{1+2\omega t}}$ solves
\begin{equation}
   d\rho_t=(3-2\omega \rho_t)\frac{dt}{1+2\omega
  t}+2\sqrt{\rho_t}\frac{dB_t}{\sqrt{1+2\omega t}}.\label{sqbes}
\end{equation}
It is easy to check that $\left(\int_0^{\frac{1}{2\omega}(e^{2\omega
  t}-1)}\frac{dB_s}{\sqrt{1+2\omega s}}\right)_t$ is a Brownian
  motion. Hence, performing a time-change in \eqref{sqbes}, one obtains that $\rho_{\frac{1}{2\omega}(e^{2\omega
  t}-1)}=e^{-2\omega t}R_{\frac{1}{2\omega}(e^{2\omega
  t}-1)}$ is a weak solution of the equation $dY_t=(3-2\omega
  Y_t)dt+2\sqrt{Y_t}\;dW_t$ satisfied by $Y_t=X_t^2$. Therefore $e^{-\omega t}\sqrt{R_{\frac{1}{2\omega}(e^{2\omega
  t}-1)}}$ is a weak solution of \eqref{eds}.\\
For $v\geq u$, $R_v$ has the same distribution as $\left(\sqrt{R_u}
  +\uu{W}^1_v -\uu{W}^1_u  \right)^2 + (\uu{W}^2_v -\uu{W}^2_u)^2 +
(\uu{W}^3_v -\uu{W}^3_u)^2$, and therefore as
  $(\sqrt{R_u}+G\sqrt{v-u})^2-2(v-u)\log(U)$ with
  $(G,U)$ independent from $R_u$. Hence for $t\geq s$, $X_{t}$ has the same
  distribution as
\begin{align*}&\left(e^{-2\omega t}\left(\left(e^{\omega
      s}X_s+\frac{G}{\sqrt{2\omega}}(e^{2\omega t}-e^{2\omega
      s})^{1/2}\right)^2-2\frac{1}{2\omega}(e^{2\omega t}-e^{2\omega
      s})\log(U)\right)\right)^{1/2}\\
&=\left(\left(e^{-\omega
      (t-s)}X_s+\frac{G}{\sqrt{2\omega}}(1-e^{-2\omega(t-s)})^{1/2}\right)^2-\frac{1}{\omega}(1-e^{-2\omega
      (t-s)})\log(U)\right)^{1/2}\end{align*}
where the couple $(G,U)$ is
      independent from $X_s$.
\subsection*{Simulation of $X_0$ with distribution $2\psi_I^2(x)1_{\{x>0\}}dx$.}
The random variable
$\frac{1}{\sqrt{2\omega}}\left(G^2-2\log(U)\right)^{1/2}$ is
distributed according to the invariant measure
$2\psi_I^2(x)1_{\{x>0\}}dx$, as suggested by letting the time increment $t-s$ tend to $+\infty$ in the previous
simulation. Indeed, $G^2-2\log(U)$ is a Gamma random
variable with density $\frac{1}{2^{3/2}\Gamma(3/2)}1_{\{z>0\}}
\sqrt{z}e^{-z/2}$. And one deduces the
density of $\frac{1}{\sqrt{2\omega}}\left(G^2-2\log(U)\right)^{1/2}$ by an easy change of variables.

\end{document}